\input ./style/arxiv-general.cfg
\documentclass[aap,MSNbibl,seceqn,dvips]{arximspdf}
\makeatletter
   \@ifpackageloaded{graphicx}{}{\usepackage{graphicx}}
\makeatother

\doi{10.1214/15-AAP1110}
\volume{26}
\issue{2}
\pubyear{2016}
\firstpage{986}
\lastpage{1028}
\docsubty{FLA}

\makeatletter
\newcommand{\rrvert}{\vert}
\newcommand{\llvert}{\vert}
\def\cal{\mathcal}

\newtheorem{theo}{Theorem}[section]

\newtheorem{prop}{Proposition}[section]
\newtheorem{lemm}{Lemma}[section]
\newproclaim{defn}{Definition}[section]
\newproclaim{assu}{Assumption}[section]

\newcommand{\Po}{{\cal P}}
\newcommand{\bH}{\mathbf{H}}
\newcommand{\cH}{{\cal H}}
\newcommand{\bQ}{\mathbf{Q}}
\newcommand{\cX}{{\cal X}}

\def\N{\mathbb{N}}
\def\R{\mathbb{R}}
\def\E{\mathbb{E}}
\def\Pr{P}

\renewcommand{\E}{\mathbb E }
\newcommand{\vold}{\pi_d} 
\newcommand{\X}{{\cal X}}
\newcommand{\A}{{\cal A}}

\newcommand{\card}{\operatorname{card}}
\newcommand{\diam}{\operatorname{diam}}
\newcommand{\st}{\mathrm{step}}
\newcommand{\K}{{\cal K}}
\newcommand{\tphi}{{\tilde{\phi}}}

\newcommand{\eps}{\varepsilon}

\def\R{\mathbb{R}}

\renewcommand{\P}{{\mathbb P}}

\def\A{{\cal A}}

\def\la{{\lambda}}

\makeatother

\begin{document}
\begin{frontmatter}

\title{Connectivity of soft random geometric graphs}
\runtitle{Soft random geometric graphs}

\begin{aug}
\author[A]{\fnms{Mathew D.}~\snm{Penrose}\corref{}\ead[label=e1]{m.d.penrose@bath.ac.uk}}
\runauthor{M.~D. Penrose}
\affiliation{University of Bath}
\address[A]{Department of Mathematical Sciences\\
University of Bath\\
Bath BA2 7AY\\
United Kingdom\\
\printead{e1}}
\end{aug}

%
\received{\smonth{4} \syear{2014}}
%
\revised{\smonth{1} \syear{2015}}

%
\begin{abstract}
Consider a graph on $n$ uniform random points
in the unit square, each pair being connected
by an edge with probability $p$ if the
inter-point distance is at most $r$.
We show that as $n \to\infty$ the probability of full connectivity is
governed by that of having no isolated vertices, itself
governed by a Poisson approximation
for the number of isolated vertices, uniformly
over all choices of $p,r$.
We determine the asymptotic probability of connectivity
for all $(p_n,r_n)$ subject to
$r_n = O( n^{-\eps})$, some $\eps>0$.
We generalize the first result to higher dimensions
and to a larger class of connection probability functions.
\end{abstract}

%
\end{frontmatter}

\section{Introduction}\label{sec1}

For certain random graph models, it
is known that the main obstacle to
connectivity is the existence of isolated vertices.
In particular, for the {Erd\H{o}s--R\'enyi }random graph $G(n,p_n)$
the probability that the graph is disconnected
but free of isolated vertices tends to
zero as $n \to\infty$, for any choice of $(p_n)_{n \geq1}$; see
\cite{ER} or \cite{Boll}, Theorem~7.3. Likewise
for the geometric graph (Gilbert graph) $G(\X_n,r_n)$ with
vertex set $\X_n$ given by a set of $n$ independently
uniformly distributed points
in $[0,1]^d$ with $d \geq2$,
and with an edge included between each
pair of vertices at distance at most $r_n$,
the probability that the graph is disconnected
but free of isolated vertices tends to
zero as $n \to\infty$, for any choice of $(r_n)_{n \in\N}$; this
follows, for example, from the results in \cite{MSTLong,kcon}.

Moreover, for both of these types of random graph (denoted $G$),
the number of isolated vertices [denoted $N_0(G)$] enjoys
a Poisson approximation for large $n$, so that with $\K$ denoting
the class of connected graphs, for large $n$ we have
%
\begin{equation}
\Pr[ G \in\K] \approx P\bigl[ N_0(G) =0 \bigr] \approx\exp(- \E
N_0). \label{eqbasic}
\end{equation}

These results have very different
proofs for geometric graphs than they do for {Erd\H{o}s--R\'enyi }graphs.
In the present paper we prove results of this kind for
a class of random graph models which generalizes both
$G(n,p)$ and $G(\X_n,r)$; we connect each pair of points
of $\X_n$ with a probability that is a function
$\phi$ of
the distance (or more generally, the displacement)
between them. The function $\phi$ is called the \emph{connection function},
and we refer to the resulting graph
as a \textit{soft} random geometric graph.\vadjust{\goodbreak}

We show that the
second approximation in (\ref{eqbasic}) holds for
soft random geometric graphs for large $n$,
uniformly over connection functions
that decay exponentially in some fixed positive power of distance,
while the first approximation in (\ref{eqbasic})
holds uniformly over connection functions
that are zero beyond a given distance, with
distance measured on the characteristic length scale of
the connection function.
For a more restricted class of connection functions, which
amount to retaining each edge of $G(\X_n,r)$ with probability $p$ in $d=2$,
we determine the limiting behavior of $\Pr[G \in\K]$
for any sequence $(r_n,p_n)_{n \geq1}$ such that there exists
$\eps>0$
with $r_n = O( n^{-\eps})$.

We also show for general $d$ that for any $(p_n)_{n \geq1}$
with $p_n \gg(\log n)/n$, if we place
the vertices of $G(n,p_n)$ at the points of
$\X_n$ and add the edges in order of
increasing Euclidean length, with high probability the
threshold for connectivity equals the threshold for
having no isolated vertices. This was previously known for
$p_n \equiv1$ \cite{kcon}.

There is substantial interest in these types of
results in the
engineering and computer science communities.
Connectivity of random geometric graphs is of interest
because of applications in wireless communications,
for example, in obtaining bounds for the
capacity of wireless networks \cite{GK,GK2}.
The ``hard'' version of the geometric
graph model (with $\phi$ the indicator of a ball centred at the origin)
is not always realistic;
communication between two nodes may not be guaranteed
even when they are close to each other
\cite{DPS,GK,MaoAnderson,YWLH}.
Also, in some cases randomness may be deliberately introduced
into the connections between nearby nodes as
a means to make the network secure \cite{KGM,Yagan,YWLH}. Among other things,
our results address a version of a conjecture of
Gupta and Kumar \cite{GK}, as
discussed at the end of Section~\ref{secmainres}.

\section{Main results}
\label{secmainres}

Throughout this paper we assume $d \in\N$ with $d \geq2$.
Given a measurable function $\phi\dvtx  \R^d \to[0,1]$
that is symmetric
[i.e., satisfies $\phi(x) =\phi(-x)$ for all $x \in\R^d$],
and given a locally finite
set $\X\subset\R^d$,
let $G_{\phi}(\X)$
be the random graph with vertex set $\X$, obtained
when each potential edge $\{x,y\}$
(with $x,y \in\X$ and $x \neq y$)
is present in the graph with probability
$\phi(x-y)$, independently of all other possible edges.

Let $\Gamma:=[0,1]^d$.
For $\la>0$ let $\cH_\la$ denote a homogeneous
Poisson point process in $\R^d$ of intensity $\lambda$,
viewed as a random subset of $\R^d$,
and let $\Po_\la: = \cH_\la\cap\Gamma$.
Given $\phi$ as above,
let $G_{\phi}(\X_n)$ and $G_{\phi}(\Po_\la)$ be the resulting
graphs as just described. We refer to $\phi$ as
the \emph{connection function}.

Soft random geometric graphs of this type
are finite-space versions of the so-called
\emph{random connection model} of continuum percolation;
for further motivation, see \cite{MR,conperc}; see also \cite{MR}, Section~1.5,
for a formal construction.

We consider various classes of connection functions
$\phi$.
Let $|\cdot|$ denote the Euclidean norm on
$\R^d$.\vadjust{\goodbreak}
Let $\Psi_d$ be the class of connection functions
$\phi$ on $\R^d$ that
satisfy
%
\begin{equation}
\phi(x) \geq\phi(y) \qquad\mbox{whenever } |x | \leq|y|. \label{0905a}
\end{equation}
%
In particular, every $\phi\in\Psi_d$ is radially symmetric,
that is, satisfies $\phi(x) =\phi(y) $ whenever $|x|=|y|$.
Condition (\ref{0905a})
is physically reasonable and is imposed
on the connection functions considered in \cite{MR}, for example.

Given a connection function $\phi$ on $\R^d$,
define the \emph{maximum value} of $\phi$ by
\[
\mu(\phi):= \sup\bigl\{\phi(x)\dvtx x \in\R^d\bigr\}.
\]
Given also $\eta>0$, let
%
\begin{equation}
\rho_\eta(\phi):= \inf\bigl\{|x|\dvtx x \in\R^d, \phi(x)
< \eta\mu(\phi) \bigr\} \label{rhodef}
\end{equation}
and also
\[
\rho_0(\phi):= \sup\bigl\{|x|\dvtx x \in\R^d, \phi(x) >
0\bigr\},
\]
which may be infinite.
%

Let $\Phi_{d,\eta}$ denote the
set of connection functions
$\phi$ on $\R^d$ such that first $\rho_\eta(\phi) \in(0,\infty
)$, second
%
\begin{equation}
\phi(x) \leq3 \eta^{-1} \mu(\phi) \exp\bigl(-\eta\bigl(|x|/
\rho_\eta(\phi)\bigr)^{\eta}\bigr),\qquad x \in\R^d
\label{eqPhi1}
\end{equation}
%
and third, $\phi\in\Psi_d$ if $d \geq3$.
Thus $\Phi_{d,\eta} \subset\Psi_d$ for $d \geq3$ but not for $d=2$.
Let $\Phi_{d,\eta}^0$ be the class of connection functions $\phi\in
\Phi_{d,\eta}$
that also satisfy
%
\begin{equation}
\rho_0(\phi) \leq\eta^{-1} \rho_\eta(\phi)
\label{eqPhi0}.
\end{equation}
%
%
For $\eta> \eta' >0$ we have
$\Phi_{d,\eta} \subset\Phi_{d,\eta'}$ and
$\Phi_{d,\eta}^0 \subset\Phi_{d,\eta'}^0$.
Condition (\ref{eqPhi1}) states that if we view $\rho_\eta(\phi)$ as
the characteristic length scale of $\phi$, then
the function $\phi(x)$ decays exponentially in the $\eta$th
power of the length of $x$, with length measured
in terms of the characteristic length scale of $\phi$.

Given $d$, define $\Psi_\st\subset\Phi_{d,1}^0 \cap\Psi_d$ by
\[
\Psi_\st:= \bigl\{\phi_{r,p}\dvtx r > 0, p \in(0,1] \bigr
\},
\]
where
for $r >0$ and $0 < p \leq1$,
we
set $\phi_{r,p}(x):= p {\mathbf1}_{[0,r]}(|x|)$.
The graph $G_{\phi_{r,p}}(\X_n)$ may be viewed
as the intersection of the (Gilbert) random geometric
graph $G(\X_n,r)$ and the {Erd\H{o}s--R\'enyi }random graph $G(n,p)$.

\emph{Rayleigh fading} functions are another class of connection functions,
where $\phi(x) = \exp(- \beta(|x|/\rho)^\gamma)$
for some fixed positive $\beta,\gamma,\rho>0$ (typically $\gamma=2$),
which is important in application; see \cite{CDG,Tse}.
Such connection functions lie in $\Phi_{d,\eta}$ for
suitable $\eta>0$, which depends on $\beta$ and $\gamma$ but
not on the length-scale $\rho$.

For any graph $G$ let
$N_0(G)$ denote the number of isolated vertices in $G$.
Also let ${\cal K}$
denote the class of connected graphs.
Our first two main results are as follows.

\begin{theo}
\label{thm2a}
Let $\eta\in(0,1]$, $k \in\N_0:=\{0,1,\ldots\}$.
Then
\[
\lim_{n \to\infty} \sup_{\phi\in\Phi_{d,\eta}}\bigl | \Pr\bigl[
N_0\bigl(G_{\phi} (\X_n)\bigr) =k \bigr] -
e^{-I_n(\phi) } I_n(\phi)^k/k!\bigr| =0,
\]
where we put
$I_n(\phi):= n \int_\Gamma\exp( - n \int_\Gamma\phi(y-x) \,dy)\,dx $.
\end{theo}

\begin{theo}
\label{thm1}
Let $\eta\in(0,1]$. Then
\begin{equation}
\lim_{n \to\infty} \sup_{\phi\in\Phi_{d,\eta}^0} \Pr\bigl[ \bigl
\{N_0\bigl(G_{\phi} (\X_n)\bigr) =0\bigr\}
\setminus \bigl\{G_{\phi}(\X_n) \in\K\bigr\}
\bigr] =0. \label{eqmain}
\end{equation}
%
\end{theo}

It is an immediate corollary of these two theorems that
for any $\eta\in(0,1]$,
\begin{equation}
\lim_{n \to\infty} \sup_{\phi\in\Phi^0_{d,\eta}} \bigl| P\bigl[
G_{\phi}(\X_n) \in\K\bigr] - \exp\bigl( - I_n(
\phi)\bigr) \bigr| =0. \label{eq2b}
\end{equation}
%
An essentially equivalent way to state the preceding results is
the following.

\begin{theo}
\label{thm2}
Let $\alpha\in[0, \infty]$
and $\eta\in(0,1]$,
and suppose $(\phi_n)_{n \in\N}$
is a sequence of connection functions
in $\Phi_{d,\eta}$,
satisfying
%
\begin{equation}
n \int_{\Gamma} \exp \biggl(- n \int_\Gamma
\phi_n(y-x) \,dy \biggr) \,dx \to\alpha \label{prcond1}
\end{equation}
as $n \to\infty$ (possibly just along some subsequence).
If $\alpha\in(0,\infty)$,
then as $n \to\infty$ (along the same subsequence if applicable),
we have for $k \in\N_0:=\{0,1,\ldots\}$ that
%
\begin{equation}
\Pr\bigl[N_0\bigl( G_{\phi_n}(\X_n)\bigr) =k
\bigr] \to e^{-\alpha} \alpha^k/k!. \label{eq1b}
\end{equation}
If $\alpha=0$, then
$\Pr[N_0( G_{\phi_n}(\X_n)) =0] \to1$,
and if
$\alpha=\infty$, then
$\Pr[N_0( G_{\phi_n}(\X_n)) =k] \to0$ for all $k \in\N_0$.
Finally, if
$\phi_n \in\Phi_{d,\eta}^0$ for all $n$,
then
%
\begin{equation}
P\bigl[G_{\phi_n}(\X_n) \in\K\bigr] \to e^{-\alpha}\qquad
\mbox{as } n \to\infty \mbox{ along the subsequence}, \label{eq1c}
\end{equation}
with
$e^{-\alpha}$ interpreted as $0$ for $\alpha= \infty$.
\end{theo}

For an example of functions that are \emph{not}
covered by our results, consider taking $\phi_n(x) = \min(1, \eps_n/|x|)$
with $\eps_n$ some sequence tending to zero. Then there is no $\eta
\in(0,1]$
such that $\phi_n \in\Phi_{d,\eta}$ for all $n$.
Another example would be if $\phi$ was the indicator of an annulus
centered at the origin; this would have $\rho_\eta(\phi)=0$ and thus not
be in $\Phi_{d,\eta}$ for any $\eta>0$.

Our definition of $\Phi_{d,\eta}$ means
we restrict attention to connection functions $\phi\in\Psi_d$
when $d \geq3$. This is because to deal with all
kinds of boundary regions of $\Gamma$ in $d \geq3$, we use the radial symmetry
of $\phi$; see Lemma~\ref{halflem}(b) below, and the result from
\cite{kcon}
or \cite{Penbk} used in its proof.
When $d=2$ the only kinds of boundary regions are either near the
corners of $\Gamma$ (a ``small'' region) or near the $1$-dimensional
edges [which can be dealt with using the condition $\phi(x) = \phi(-x)$;
see Lemma~\ref{halflem}(a) below],
so we do not require $\phi\in\Psi_2$ for the results above.


Given $r \geq0$ and $ p \in(0,1]$ and finite $\X\subset\Gamma$,
write $G_{r,p}(\X)$ for $G_{\phi_{r,p}}(\X)$.
Given $p$,
a natural coupling of all the graphs $G_{r,p}(\X_n), r \geq0$,
goes as follows:
let $G_{r,p}(\X_n)$ be the subgraph of
$G_{\sqrt{d},p}(\X_n)$,
with vertex set $\X_n$,
and edge set consisting of
all edges of Euclidean length at most
$r$. With this coupling,
$G_{r,p}(\X_n)$ is a subgraph of
$G_{s,p}(\X_n)$ whenever $r \leq s \leq\sqrt{d}$.
Given $p$, define the thresholds
$
\tau_n(p):= \inf\{r\dvtx G_{r,p}(\X_n) \in\K\}$,
and
$
\sigma_n(p):= \inf\{r\dvtx N_0( G_{r,p}(\X_n)) = 0 \}$,
with the infimum of the empty set interpreted as $+\infty$.
Clearly $\sigma_n(p) \leq\tau_n(p)$ almost surely.
Our next result gives an
asymptotic equivalence of these two thresholds.

\begin{theo}
\label{threshthm}
Given any $[0,1]$-valued sequence
$(p_n)_{n \in\N}$ with $n p_n/\break \log n \to\infty$
as $n \to\infty$, it is the case that
\[
\lim_{n \to\infty} P\bigl[\tau_n(p_n)
= \sigma_n(p_n)\bigr] =1.
\]
\end{theo}

In the case where $d =2$ and $\phi_n \in\Psi_2 \cap\Phi_{2,\eta}$
for some $\eta\in(0,1]$,
we shall make Theorem~\ref{thm2} more explicit,
by characterizing those sequences $\phi_n$ which satisfy
(\ref{prcond1}).
Setting $p_n:= \mu(\phi_n)$,
we find that the
main contribution to the integral in (\ref{prcond1})
comes from $x$ in the interior of $\Gamma$ when $p_n \gg(1/\log n)$,
while the main contribution comes from $x$ near the boundary
but not the corners of $\Gamma$ when $ n^{-1/3} (\log n)^{-1} \ll p_n
\ll1/\log n$,
and the main contribution comes from $x$ near the corners of $\Gamma$
when $p_n \ll n^{-1/3} (\log n)^{-1}$.

We state this more precisely in Theorem~\ref{thconds} below, which requires further notation.
Given real-valued functions $f,g$,
recall that
$f(n) = \omega(g(n))$ means
$g(n) = o(f(n))$ (as $n \to\infty$),
$f(n) = \Omega(g(n))$ means
$g(n) = O(f(n))$ and
$f(n) = \Theta(g(n))$ means
$f(n) = O(g(n))$ and $g(n) = O(f(n))$.
Finally $f(n) \sim g(n)$ means $f(n) = (1+o(1)) g(n)$.
For any connection function $\phi$ we set
%
\begin{equation}
I(\phi): = \int_{\R^d} \phi(x) \,dx. \label{Idef}
\end{equation}
If $\eta\in(0,1]$ and $\phi\in\Phi_{2,\eta}$, then set
%
\begin{eqnarray}
J_1(\phi)&:=& J_1(\phi,\eta):= \mu(\phi)^{-1}
\int_{0}^\infty\phi\bigl( \bigl(
\rho_\eta(\phi) t,0\bigr) \bigr) \,dt; \label{J1def}
\\
J_2(\phi)&:=& J_2(\phi,\eta):= \mu(\phi)^{-1}
\int_{0}^\infty\phi\bigl( \bigl(
\rho_\eta(\phi) t,0\bigr) \bigr) 2\pi t \,dt. \label{J2def}
\end{eqnarray}
For $\eta\in(0,1]$
and $\phi\in\Psi_2 \cap\Phi_{2,\eta}$, we have
$I(\phi) = \mu(\phi) \rho_\eta(\phi)^2 J_2(\phi)$,
and for
$\phi\in\Psi_\st$ we have
$J_1(\phi) =1$ and
$J_2(\phi) =\pi$.

The integrals $J_1(\phi)$ and $J_2(\phi)$ may be viewed
as measures of the ``shape'' of $\phi$, separate from $\mu(\phi)$
and $\rho_\eta(\phi)$, which
measure the vertical and horizontal ``scale'' of $\phi$, respectively.
Note that for $\eta\in(0,1]$
and $i=1,2$, we have
%
\begin{equation}
0 < \inf_{\phi\in\Psi_2 \cap\Phi_{2,\eta}} J_i(\phi,\eta) \leq\sup
_{\phi\in\Psi_2 \cap\Phi_{2,\eta}} J_i(\phi,\eta) < \infty. \label{Jbds}
\end{equation}

\begin{theo}
\label{thconds}
Let $\eta\in(0,1]$, $\alpha\in(0,\infty)$.
Suppose $d =2$ and
$\phi_n \in\Phi_{2,\eta} \cap\Psi_2$ for $n \in\N$.
Set $r_n:= r_\eta(\phi_n)$ and
$p_n:= \mu(\phi_n)$.
Then (\ref{prcond1}) holds under any of the following conditions
as $n \to\infty$:
\begin{longlist}[(1)]
\item[(1)] $p_n = \omega(1/\log n)$ and
$
n I(\phi_n) - \log n \to- \log\alpha$;
\item[(2)]
$p_n = o(1/\log n)$ and $p_n = \omega( n^{-1/3} (\log n)^{-1})$
and
%
\begin{equation}
n I(\phi_n) = \log \biggl( \frac{4 J_2(\phi_n)
}{\alpha^2J_1(\phi_n)^2} \biggr) + \log \biggl(
\frac{n}{p_n} \biggr) - 
\log\log \biggl( \frac{n}{p_n} \biggr)
+ o(1); \label{1018a}
\end{equation}
\item[(3)]
$p_n = o(n^{-1/3} (\log n)^{-1}) $ and $r_n = n^{-\Omega(1)}$ and
\[
n I(\phi_n) = 4 \bigl( \log(1/p_n) 
- \log
\log(1/p_n) + \log\bigl(J_2(\phi_n)/ \bigl(
\alpha J_1(\phi_n)^2\bigr) \bigr) \bigr)
+o(1).
\]
\end{longlist}
\end{theo}

We also deal with the boundary cases
$p_n = \Theta(1/\log n)$ and $p_n =  \break \Theta(
n^{-1/3} (\log n)^{-1}) $; see Theorems
\ref{lemrp5} and \ref{lemrp4}.

We now discuss related work and open problems.
Note that (\ref{eq1b}) [but not (\ref{eq1c})] of Theorem~\ref{thm2} was
already proved
by Yi et al. \cite{YWLH} in the special case
with $d=2$ and $\phi_n \in\Psi_\st$
under the condition $p_n = \omega(1/\log n)$.
Here we are considering a much more general class
of sequences of connection functions $\phi_n$.

For a discussion of these problems from a statistical
physics viewpoint via formal series expansions
and for further discussion of motivation, see
Coon et al. \cite{CDG}.
The methods of
Krishnan et al.
\cite{KGM}
(see Remark~3 of that paper) could be used
to give some limiting inequalities for the probability
of connectivity
in the special case of connection functions in
$\Psi_\st$ [whereas our
(\ref{eq2b}) provides a limiting \emph{equality} for a more general
class of connection functions].
The main concern in \cite{KGM} is with a certain nonindependent randomization
(random key graphs) to determine which of the edges (below the
threshold radius)
are present, which is of interest
from an engineering perspective; see also \cite{Yagan}.
It would be interesting to try to
extend our results to these random key graphs.

A related random graph model is the
\emph{bluetooth graph};
this is a subgraph of the ``hard'' random geometric graph
with edges selected at random according to
a restriction on vertex degrees. See \cite{BDFL} for
results on connectivity of bluetooth graphs.

Another related problem is that of \emph{Hamiltonicity}.
Analogously to (\ref{eqmain}),
one might speculate that
for large $n$,
the probability that
$G_{\phi}(\X_n)$ is non-Hamiltonian while having
minimum degree at least 2, might vanish
uniformly over connection functions in $\Psi_\st$
(or indeed, connection functions in $\Phi_{d,\eta}^0$).
For the more restricted class of connection functions
of ``hard'' random geometric graphs, this was proved in
\cite{BBKMW}.
Some of the ideas of proof in the present paper are related
to methods used in \cite{BDFL} and in \cite{BBKMW}.

Given $k \in\N$, and given a graph $G$,
let $N_{<k}(G)$ be the number of vertices of
$G$ of degree less than $k$, and let
$\K_k$ be the class of $k$-connected graphs.
In view of the results from \cite{kcon}, one might expect
(\ref{eqmain})
to hold with
$N_0$ replaced by $N_{<k}$ and
$\K$ replaced by $\K_k$, for any fixed $k \in\N$.

In a much-cited paper, Gupta and Kumar \cite{GK} conjectured that
if $d=2$, $\X_n$ consists of $n$
points uniformly distributed in a disk of unit area
(rather than the unit square considered here),
and $\phi_n = \phi_{r_n,p_n}$,
then $\P[ G_{\phi_n}(\X_n) \in\K] \to1$
if and only if $n \pi r_n^2 p_n -\log n \to\infty$.
Our results
(Theorems \ref{thm2}, \ref{thconds} and
\ref{lemrp5})
address the corresponding
conjecture for points in the unit square, showing
that under the additional assumption that
$p_n = \Omega(1/\log n) $,
the conjecture is true and also
$\P[ G_{\phi_n}(\X_n) \in\K] \to0 $ if
$n \pi r_n^2 p_n -\log n \to- \infty$.
Our results also show that if
$p_n = \omega(1/\log n) $
and
if $n \pi r_n^2 p_n -\log n \to\beta\in\R$, then
$\P[ G_{\phi_n}(\X_n) \in\K] \to\exp(-e^{-\beta})$.

However, if one assumes instead that $p_n = o(1/\log n)$
and $p_n = \break  \omega(n^{-1/3} (\log n)^{-1})$
and (\ref{1018a}) holds, then it is easily
verified that
$n \pi r_n^2 p_n -\log n \to\infty$,
but our results show that
$\P[ G_{\phi_n}(\X_n) \in\K] $ tends
to a limit strictly between 0 and 1, so the conjecture
fails. Essentially, this is because, in this
case, the mean number of isolated vertices in
the interior of $\Gamma$ tends to zero, but
the mean number of isolated vertices near the
boundary does not. In this regime the corner
effects are not the most important, and
we would expect something similar to hold
in the unit disk, as considered in \cite{GK}.
More generally, it would be of interest to extend our
results to the case of other shaped regions such as
smoothly bounded regions, but this would be a nontrivial task because
the boundary effects
can be quite strong [essentially because of
the exponential factor in the expression on the left of
(\ref{prcond1})].

The remaining sections of the paper are organized as follows.
In Section~\ref{Poissec} we prove
Theorem~\ref{Poissecthm}, which is
a Poissonized version
of Theorem~\ref{thm2a}
(i.e., one with the point process $\X_n$ replaced by $\Po_n$),
of interest in its own right.
In Sections \ref{Sec:psmall} and \ref{Sec:plarge},
we prove Theorem~\ref{prop0711c}, which is (loosely speaking)
a Poissonized version of
Theorem~\ref{thm1}, also of interest in its own right.

In Section~\ref{secdepo}, we shall de-Poissonize,
thereby completing the proof of Theorems \ref{thm2a}, \ref{thm1}
and \ref{thm2}.
In Section~\ref{secequiv} we prove Theorem~\ref{threshthm}.
In Section~\ref{secchoice},
we prove
Theorems \ref{thconds},
\ref{lemrp5} and \ref{lemrp4}.

We conclude this section with some remarks on the proofs.
As we have mentioned, many of the results presented here
might naturally be conjectured in view of known
results for random ``hard'' geometric graphs
\cite{MSTLong,kcon}, for {Erd\H{o}s--R\'enyi }random graphs
\cite{ER,Boll} and a (slightly weaker)
explicit conjecture along these lines given in \cite{GK}.
These references date back to the last century,
but the conjectures have not been proved before now,
despite the considerable influence of
Gupta and Kumar \cite{GK}
in the applied literature; see, for example, the discussion in \cite{Yagan}.

We believe that there are two reasons for this. One is
that different arguments are used to prove
these results depending on whether or not
$\mu(\phi_n)$ tends to zero faster than a certain rate.
The division between
Sections \ref{Sec:psmall} and \ref{Sec:plarge} reflects this,
and Section~\ref{Poissec} is also divided along these lines.
The balance between geometrical and combinatorial
arguments is different in these different settings.

The other reason is that the proof is not just a matter
of reassembling known arguments.
For example, a part of the argument
is concerned with ruling out the possibility that there
are two large disjoint components.
For ``hard'' geometric graphs
\cite{MSTLong,kcon}, any two such components
are separated by a connected region of empty
space, and one can use discretization, spatial
independence and path-counting arguments
directly. In the present ``soft'' case, however, the physical
separation of components is not at all obvious.
Instead, we proceed more indirectly via
a notion of local good behavior of our point
process (the ``blue cubes'' of
Section~\ref{Subsec:LC}) with finite-range dependence,
after which we can use path-counting arguments to
establish that there is a single giant region of ``blue cubes''
corresponding to a single large component of
our graph.


\section{Poisson approximation}
\label{Poissec}
In this section we prove
the following Poissonized version of
Theorem~\ref{thm2a} (we shall de-Poissonize in Section~\ref{secdepo}).

\begin{theo}
\label{Poissecthm}
Let $\alpha>0$ and $\eta\in(0,1]$.
Suppose $(\la(n))_{n \in\N}$ is
an increasing $(0,\infty)$-valued sequence
that tends to $\infty$ as $n \to\infty$,
and $(\phi_\la)_{\la>0}$ is a collection of
connection functions in $\Phi_{d,\eta}$.
Suppose that
as $\la\to\infty$ along the sequence $(\la(n))$,
we have
%
\begin{equation}
\la\int_{\Gamma} \exp \biggl( - \la\int_{\Gamma}
\phi_\la(y-x) \,dy \biggr) \,dx \to\alpha. \label{prcond0}
\end{equation}
Then for $k \in\N_0$ we have
as $\la\to\infty$
along the same sequence,
that
%
\begin{equation}
\Pr\bigl[N_0\bigl( G_{\phi_\la}(\Po_\la)\bigr) =k
\bigr] \to e^{-\alpha} \alpha^k/k!. \label{eq1}
\end{equation}
\end{theo}

Our strategy of proof is as follows.
When
$p_\lambda:= \mu(\phi_\lambda)$ is ``small,''
we use the method of moments, the Mecke formula
(\ref{Mecke}) and Bonferroni bounds. When $p_\lambda$
is ``big'' we shall proceed by
the Chen--Stein method for Poisson approximation
of $N_0(G_{\phi_\lambda}(\Po_\lambda))$,
which may be approximated
(via discretization of space)
by a sum of ``mostly independent''
indicator functions.

In proving (\ref{eq1}), we shall use the following notation.
We write \textit{with high probability} or \textit{w.h.p.}
to mean \textit{with probability tending to \textup{1} as $\la\to\infty$}.
All asymptotic statements are taken to be as
$\la\to\infty$ along the sequence $\la(n)$ mentioned
in Theorem~\ref{Poissecthm}.
Also, for $A,B \subset\R^d$ we write
$A \oplus B$ for $\{x+y\dvtx x \in A,y \in B\}$ (Minkowski
addition of sets).

For any finite (deterministic) $\A\subset\R^d$, and any $\phi\in
\Phi_{d,\eta}$,
set
%
\begin{equation}
h_{\phi}(\A):= P\bigl[ G_{\phi}(\A) \in\K\bigr],
\label{hphidef}
\end{equation}
and for any
$y \in\R^d$ with $y \notin\A$,
set
%
\begin{eqnarray}\label{gphidef}
g_{\phi}(y,\A)&:=& 1 - \prod_{x \in\A} \bigl(1-
\phi(y-x)\bigr)
\nonumber
\\[-8pt]
\\[-8pt]
\nonumber
&=& P\bigl[y \mbox{ is nonisolated in } G_{\phi}\bigl(\A
\cup\{ y\}\bigr)\bigr].
\end{eqnarray}

The left-hand side of (\ref{prcond0}) equals $\E N_0(G_{\phi_\la})$.
This is a consequence of the following formula, which
we shall use repeatedly.
Suppose $k \in\N$
and $f$ is a
measurable nonnegative function defined on
$(\R^d)^k \times{\cal G}_k$ where
${\cal G}_k$ is the space of all graphs on vertex set
$\{1,\ldots,k\}$. Then given a connection
function $\phi$, for $\la>0$ we have
%
\begin{eqnarray}\label{Mecke}
&&\E\sum^{\neq}_{X_1,\ldots,X_k \in\Po_\la} f\bigl(X_1,
\ldots,X_k,G_\phi(\Po_\la)|_{X_1,\ldots,X_k}
\bigr) {\mathbf1}_{{\mathbf D}_\phi(X_1,\ldots,X_k;\Po_\la) }\nonumber\\
&&\qquad= \la^k \int_\Gamma
\,dx_1 \cdots \int_\Gamma \,dx_k \E\bigl[f\bigl(x_1,\ldots,x_k,G_\phi
\bigl(\{x_1,\ldots,x_k\}\bigr)\bigr)\bigr] \\
&&\qquad\quad{}\times\exp \biggl( -
\la\int_\Gamma g_\phi\bigl(y; \{x_1,
\ldots,x_k\} \bigr) \,dy \biggr), \nonumber
\end{eqnarray}
where the sum is over all ordered $k$-tuples of distinct
points of
$\Po_\la$, and
$G_\phi(\Po_\la)|_{X_1,\ldots,X_k} $ is the subgraph
of $G_\phi(\Po_\la)$ induced by vertex set $\{X_1,\ldots,X_k\}$
with the vertex $X_i$ given the label $i$ for each $i$,
and
${\mathbf D}_\phi(X_1,\ldots,X_k;\Po_\la) $ is the event that
there is no edge of $G_\phi(\Po_\la)$ between any vertex in $\{
X_1,\ldots,X_k\}$
and any vertex in $\Po_\la\setminus\{X_1,\ldots,X_k\}$.

Formula (\ref{Mecke}) is related to
the Slivnyak--Mecke formula in the theory of Poisson
processes; here we just call it the \emph{Mecke formula}.
It can be proved by conditioning on the number of points
of $\Po_\la$; see the proofs of \cite{Penbk}, Theorem~1.6 and
\cite{conperc}, Proposition~1.

We shall use the following inequality more than once.
Given connection function $\phi$ and given
$x,x_1,\ldots,x_k \in\Gamma$, by the Bonferroni bound
\[
g_{\phi} \bigl(x;\{x_1,\ldots,x_k\}
\bigr) \geq 
 \Biggl( \sum_{i=1}^k
\phi(x-x_i) \,dx \Biggr) - 
\sum
_{1\leq i < j \leq k} \phi(x-x_i) \phi(x-x_j),
\]
so integrating over $x \in\Gamma$, we obtain
%
\begin{equation}
\int_\Gamma g_{\phi} \bigl(x;
\{x_1,\ldots,x_k\} \bigr) \,dx \geq 
 \Biggl( \sum
_{i=1}^k \int_\Gamma
\phi(x-x_i) \,dx \Biggr) 
- k^2 \mu(\phi) I(\phi).
\label{0705a}
\end{equation}
%

Let $\bH$ denote the half-space $[0,\infty) \times\R^{d-1}$,
and let $\bQ$ denote the
orthant $[0,\infty)^d$.
For $x \in\bQ$ let $\bQ_x: = \{y \in\bQ\dvtx  \|x\|_1 \leq\|y\|_1\}$,
where $\|\cdot\|_1$ is the $\ell_1$ norm.

\begin{lemm}
\label{halflem}
Let $\eta\in(0,1]$
and $\phi\in\Phi_{d,\eta}$. Then:
\textup{(a)} if $d=2$, for any $x= (x_1,x_2) \in\bH$ and
$y= (y_1,y_2) \in\bH$ with $x_1 \leq y_1$,
and $r \in[\rho_\eta(\phi),\infty]$,
setting $\phi^{(r)} (x):= \phi(x) {\mathbf1}_{[0,r]}(|x|)$
we have
\[
\int_{\bH} \bigl( g_{\phi^{(r)}}\bigl(z,\{x,y\}\bigr) -
\phi^{(r)}(z-x) \bigr) \,dz \geq ( \eta/4) \mu(\phi) \rho_\eta(
\phi) \min\bigl(|y-x|,\rho_\eta(\phi) \bigr);
\]
\textup{(b)} if $d \geq3$, and $x \in\bQ, y \in\bQ_x$,
then
%
\begin{eqnarray}\label{150121a}
&&\int_{\bQ} \bigl( g_\phi\bigl(z,\{x,y\}\bigr) -
\phi(z-x) \bigr) \,dz
\nonumber
\\[-8pt]
\\[-8pt]
\nonumber
&&\qquad\geq \eta_1 \mu(\phi) \rho_\eta(
\phi)^{d-1} \min\bigl(|y-x|,\rho_\eta(\phi) \bigr),
\end{eqnarray}
where $\eta_1 >0$ is a constant depending only on $d$ and $\eta$.
\end{lemm}

\begin{pf}
(a)
Let us assume $x_2 \leq y_2$ (the other case may be treated
similarly).
For any $z \in\R^2$,
since
$ g_{\phi^{(r)}}(z,\{x,y\}) - \phi^{(r)} (z-x) = (1- \phi^{(r)}(z-x)
) \phi^{(r)}
(z-y) $,
we have
$g_{\phi^{(r)}}(z,\{x,y\}) - \phi^{(r)}(z-x)
\geq(\phi^{(r)}(z-y) - \phi^{(r)}(z-x))_+$.
Therefore it suffices to prove
%
\begin{eqnarray}\label{150121b}
&&\int_{\bH} \bigl( \phi^{(r)}(z-y) -
\phi^{(r)}(z-x) \bigr)_+ \,dz
\nonumber
\\[-8pt]
\\[-8pt]
\nonumber
&&\qquad \geq ( \eta/4) \mu(\phi)
\rho_\eta(\phi) \min\bigl(|y-x|,\rho_\eta(\phi) \bigr).
\end{eqnarray}
Now
\begin{eqnarray*}
\int_{\bH} \bigl( \phi^{(r)}(z-y) -
\phi^{(r)}(z-x) \bigr)_+ \,dz &\geq&\int_{\{y\} \oplus\bQ} \bigl(
\phi^{(r)}(z-y) - \phi^{(r)}(z-x) \bigr) \,dz
\\
&=& \int_{\bQ} \phi^{(r)}(w) \,dw - \int
_{\{y-x\} \oplus\bQ} \phi ^{(r)}(w) \,dw
\\
&=& \int_{\bQ\setminus( \{y-x\} \oplus\bQ) } \phi^{(r)}(w) \,dw.
\end{eqnarray*}
If $|y-x| \leq\rho_\eta(\phi)$, then the region
$\bQ\setminus( \{y-x\} \oplus\bQ) $
contains either the rectangle
$[0,|y-x|/2] \times[0,\rho_{\eta}(\phi)/2]$
or the rectangle
$ [0,\rho_\eta(\phi)/2] \times[0,|y-x|/2] $ (or both), and the
function $\phi^{(r)}$ exceeds $\eta\mu(\phi)$ on either of these
rectangles, so that
$\int_{\bQ\setminus( \{y-x\} \oplus\bQ) } \phi^{(r)}(w) \,dw \geq
\eta|y-x| \rho_\eta(\phi) \mu(\phi)/4$.\vspace*{1pt}

If $|y-x| \geq\rho_\eta(\phi)$, then the region
$\bQ\setminus( \{y-x\} \oplus\bQ) $
contains the square $[0, \rho_\eta(\phi)/2]^2$, so that
$\int_{\bQ\setminus( \{y-x\} \oplus\bQ) } \phi^{(r)}(w) \,dw \geq
\eta\rho_\eta(\phi)^2 \mu(\phi)/4$.
This gives us~(\ref{150121b}).

(b)
Now suppose $d \geq3$ (so $\phi\in\Psi_d$ by definition of $\Phi
_{d,\eta}$).
For $x,y \in\bQ$, we have
by Fubini's theorem and
(\ref{rhodef}) that
%
\begin{eqnarray} \label{fromFub2}
\int_\bQ\bigl( g_{\phi}\bigl(z,\{x,y\}\bigr) -
\phi(z-x) \bigr) \,dz& =& \int_0^1 \int
_\bQ( {\mathbf1}_{\{ g_{\phi} (z,\{x,y\}) \geq t \}} - {\mathbf1}_{\{ \phi(z-x) \geq t \}}
) \,dz \,dt
\nonumber
\\
&\geq&\int_{0}^{\eta\mu(\phi)} \int
_\bQ( {\mathbf1}_{\{ \phi(z-y) \geq t \}} - {\mathbf1}_{\{ \phi(z-x) \geq t \}}
)_+ \,dz \,dt
\\
&= &\int_{0}^{\eta} \bigl|\bQ\cap B\bigl(y;
\rho_{u}(\phi)\bigr) \setminus B\bigl(x;\rho_{u} (\phi)
\bigr) \bigr| \mu(\phi) \,du,\nonumber
\end{eqnarray}
where $|\cdot|$ denotes Lebesgue measure or
the Euclidean norm according to context.

For $ u \leq\eta$,
we have $
\rho_{u} (\phi) \geq\rho_\eta(\phi)$.
Also, there is
a constant $\eta_2 >0$ (dependent on $\eta$ and $d$) such that
$ |\bQ\cap B(y;1) \setminus B(x;1)| \geq\eta_2
\min( |y-x|,1)$
for any $x,y \in\bQ$ with
$\|x\|_1 \leq\|y\|_1$;
see \cite{Penbk}, Proposition~5.16
or \cite{kcon}, Proposition~2.2.
Hence for $x \in\bQ$, $y \in\bQ_x$
and $u \in(0,\eta]$,
by scaling
\begin{eqnarray*}
\bigl| \bQ\cap B\bigl(y;\rho_{u} (\phi)\bigr) \setminus B\bigl(x;
\rho_{u} (\phi)\bigr) \bigr|& \geq&\bigl(\rho_{u}(\phi)
\bigr)^d \eta_2 \min \biggl( \frac{|y-x|}{\rho_{u}(\phi) }, 1 \biggr)
\\
&\geq&\eta_2 \rho_\eta(\phi)^{d-1} \min
\bigl(|y-x|,\rho_\eta(\phi)\bigr).
\end{eqnarray*}
%
Putting this into (\ref{fromFub2}) gives us result (\ref{150121a})
with $\eta_1 = \eta_2 \eta$.
\end{pf}

Given $\eta\in(0,1]$ and given $(\phi_\la)_{\la>0}$
with each $\phi_\la\in\Phi_{d,\eta}$, for $\la>0$ we set
%
\begin{equation}
p_\la:= \mu(\phi_\la);\qquad r_\la:=
\rho_\eta(\phi_\la). \label{plrldef}
\end{equation}
Recall from (\ref{Idef}) that $I(\phi): = \int_{\R^d} \phi(x) \,dx$
for any connection function $\phi$. Without loss of generality for
the purpose of proving Theorem~\ref{Poissecthm},\vspace*{1pt} we can
and do assume for all $\la$ that
$\rho_0(\phi_\la) \leq\sqrt{d}$, so that also
$r_\la\leq\sqrt{d}$. Note that if (\ref{prcond0}) holds, then
%
\begin{equation}
\la I(\phi_\la) = \Theta ( \log\la ), \label{Thetacond}
\end{equation}
and therefore by (\ref{plrldef}),
%
\begin{equation}
\la p_\la r_\la^d = \Theta(\log\la).
\label{fromThetacond}
\end{equation}

Theorem~\ref{Poissecthm} follows from the next two lemmas,
dealing separately with the case with $p_\la= o(1/\log\la)$ and
the case with $p_\la= \omega(1/(\log\la)^2)$.
In the first case, we use the method of moments.
For $m, r \in\N$ we write $(m)_r$ for the descending factorial
$m(m-1)\cdots(m-r+1)$.

\begin{lemm}
\label{poislem2}
Let $\alpha\in(0,\infty)$, $\eta\in(0,1]$.
Suppose $\phi_\la\in\Phi_{d,\eta}$ for all $\la$ and
$(\phi_\la)_{\la>0}$ satisfy (\ref{prcond0}), and
that $p_\la= o(1/\log\la)$.
Then
(\ref{eq1}) holds.
\end{lemm}

\begin{pf}
Set $N_0:= N_0(G_{\phi_\la}(\Po_\la))$.
Let $k \in\N$.
For finite $A \subset\R^d$, let $u_\lambda({\cal A})$ denote the probability
that $G_{\phi_\lambda}({\cal A})$ has no edges.
By
the Mecke formula (\ref{Mecke}),
\begin{eqnarray*}
\E\bigl[(N_0)_k \bigr] &=& \la^k 
\int\cdots\int u_\la \bigl(\{x_1,\ldots,x_k\}
\bigr) \\
&&\hspace*{43pt}{}\times\exp \biggl(- \la\int g_{\phi_\la}\bigl(x, \{x_1,
\ldots,x_k\} \bigr) \,dx \biggr)
\,dx_1 \cdots \,dx_k ,
\end{eqnarray*}
where all integrals
are over $\Gamma$, unless specified otherwise.
By the union bound,
$u_\la(\{x_1,\ldots,x_k\}) \geq1 - {k \choose2} p_\la$,
and also
$g_{\phi_\la}(x, \{x_1, \ldots,x_k\} ) \leq
\sum_{i=1}^k \phi_{\la}(x-x_i)$.
Hence
%
\begin{eqnarray}\label{0711c}
&&\E\bigl[(N_0)_k \bigr]\nonumber\\
&&\qquad\ge 
\bigl(1 -
k^2 p_\la\bigr) \la^k 
\int\cdots
\int \exp \Biggl(- \la\int\sum_{i=1}^k
\phi_{\la}(x- x_i ) \,dx \Biggr) \,dx_1 \cdots
\,dx_k
\\
&&\qquad= \bigl(1+o(1)\bigr) (\E N_0)^k. \nonumber
\end{eqnarray}
Also, by (\ref{0705a}),
we
have
%
\begin{eqnarray}\label{0711b}
&&\E\bigl[(N_0)_k \bigr]\nonumber\\
&&\qquad\leq \la^k
\int\cdots\int \exp \Biggl( \la k^2 p_\la I(
\phi_\la) - \la\int\sum_{i=1}^k
\phi_{\la}(x - x_i ) \,dx \Biggr) \,dx_1 \cdots
\,dx_k
\\
&&\qquad= \bigl(1+o(1)\bigr) (\E N_0)^k, \nonumber
\end{eqnarray}
where the last line is due to the fact that $\la p_\la I(\phi_\la)
= O( p_\la\log\la) \to0$,
by (\ref{Thetacond})
and our assumption on $p_\la$.

By (\ref{0711c}), (\ref{0711b}) and assumption
(\ref{prcond0}), we have that $\E[ (N_0)_k] \to\alpha^k$,
and therefore by the method of moments
(see, e.g., Theorem~1.22 of \cite{Boll}),
we have
Poisson convergence
(\ref{eq1}).
\end{pf}

For the second case with
$p_\la= \omega((\log\la)^{-2})$,
we use the Poisson approximation method from
\cite{MSTLong}.
This method has the potential to provide error bounds,
but this is not our main focus here.
For $x \in\R^d$ and $r>0$ set $B(x;r)$ to be the ball
$\{y \in\R^d\dvtx |x-y| \leq r \}$.
Given $\eta\in(0,1]$, set
\[
K(\eta): = \int_{\R^d} 3 \eta^{-1} \exp\bigl(-
\eta|x|^{\eta} \bigr) \,dx. 
\]
Note that
$K(1) \leq K(\eta) < \infty$, and $K(1) = 6\pi$ if $d=2$, and that
by (\ref{eqPhi1}) and~(\ref{Idef}),
\begin{equation}
I(\phi) \leq\mu(\phi) \bigl(\rho_\eta(\phi)\bigr)^d K(
\eta), \qquad\phi\in\Phi_{d,\eta}. \label{0729b}
\end{equation}

\begin{lemm}
\label{poislem1}
Suppose for some $\eta\in(0,1]$ and
$\alpha\in(0,\infty)$ that
$\phi_\la\in\Phi_{d,\eta}$ for all $\la>0$ and $\phi_\la$
satisfy (\ref{prcond0}). Suppose
$p_\la= \omega(1/(\log\la)^2)$.
Then (\ref{eq1}) holds.
\end{lemm}

\begin{pf}
Assume $r_\la\leq\sqrt{d}$.
It follows from (\ref{prcond0}) that
(\ref{Thetacond}) and (\ref{fromThetacond}) hold.
Hence by our condition on $p_\la$ we have
%
\begin{equation}
r_\la^d = \Theta\bigl((\log\la) /(\la p_\la)
\bigr) = o \bigl( (\log\la)^{3}\la^{-1}\bigr). \label{0729a}
\end{equation}
By (\ref{fromThetacond}), we can (and do) choose
$\delta> 0$ with $\lambda p_\lambda r_\lambda^d > \delta\log\lambda$
for all $\lambda$.
Let $\eps>0$ be fixed with
$\eps< \eta/(4K(\eta))$ if $d =2$, and with
$\eps< \min( 2^{-d} \pi_d \eta/K(\eta), \eta_1 \delta) $
if $d \geq3$, where $\eta_1$ is as in Lemma~\ref{halflem}(b).
Truncate $\phi_\la$ by setting $\tilde{\phi}_\la(x):=
\phi_\la(x) {\mathbf1}_{[0,r_\la^{1-\eps}]} (|x|)$ for $x \in\R^d$.
Couple
$G_{\phi_\la}(\Po_\la)$ and
$G_{\tilde{\phi}_\la}(\Po_\la)$
in the following natural way:
starting with
$G_{\phi_\la(\Po_\la)}$, remove all edges
of Euclidean length greater than $r_\la^{1-\eps}$ to obtain
$G_{\tilde{\phi}_\la}(\Po_\la)$.

We claim next that (\ref{prcond0}) holds with $\phi_\la$
replaced by $\tilde{\phi}_\la$, that is,
%
\begin{equation}
\la\int_{\Gamma} \exp \biggl( - \la\int_{\Gamma}
\tilde{\phi}_\la(y-x) \,dy \biggr) \,dx \to \alpha. \label{trunc}
\end{equation}
Indeed, by the Mecke formula (\ref{Mecke}) the absolute value of the
difference between the left-hand side of (\ref{trunc})
and that of (\ref{prcond0}) is bounded by
the mean number of vertices having at least one
incident edge in $G_{\phi_\la}(\Po_\la)$ of length
at least $r_\la^{1-\eps}$,
and hence by twice the expected number of such edges.
However,
by (\ref{eqPhi1})
the expected number
of such edges is
$O(\la^2 \exp(-\eta r_\la^{- \eps\eta}))$,
which is $O(\la^2 \exp(-\eta\la^{\eps\eta/(2d)}))$
by~(\ref{0729a}), and therefore tends to zero.

Let $\Gamma'_\la$ be the set of $x \in\Gamma$ distant more than
$4r_\la^{1-\eps}$ in the $\ell_\infty$ norm from the corners of
$\Gamma$.
Let $\tilde{N}_0(\la)$ be the number of isolated vertices
of $G_{\tilde{\phi}_\la}(\Po_\la)$ that are located in $\Gamma
'_\la$.
Then
we claim that
%
\begin{equation}
\E\bigl[ \bigl| N_0\bigl(G_{\phi_\la}(\Po_\la)\bigr)-
\tilde{N}_0 (\la) \bigr|\bigr] \to0\qquad \mbox{as } \la \to\infty. \label{01claim}
\end{equation}
%
To see this, observe
first that
$\E[ | N_0(G_{\phi_\la}(\Po_\la)) - N_0(G_{\tilde{\phi}_\la
}(\Po_\la)) |] $
is bounded by twice the expected number of edges in
$G_{\phi_\la}(\Po_\la)$ of Euclidean length greater than
$r_\la^{1-\eps}$, which
tends to zero as discussed above.
Second, observe that
for all $x \in\Gamma$, by (\ref{0729b}) we have
\[
\int_\Gamma\tilde{\phi}_\la(y-x) \,dy
\geq2^{-d} \vold r_\la^d \eta p_\la
\geq I(\phi_\la) 2^{-d} \vold\eta/K(\eta),
\]
and $e^{-\la I ( \phi_\la)} = O(1/\la)$ by (\ref{prcond0}),
so that $\exp(-\la\int_{\Gamma}\tilde{\phi}_\la(y-x) \,dy) = \break  O(
\la^{-
2^{-d} \vold
\eta/K(\eta)})$, uniformly over $x \in\Gamma$.
Hence the expected number of isolated vertices
of $G_{\tilde{\phi}_\la}(\Po_\la)$ lying in
$\Gamma\setminus\Gamma'_\la$ is
$
O( r_\la^{d(1-\eps)} \la^{1 - 2^{-d } \vold\eta/K(\eta)})
$
which tends to zero by (\ref{0729a}).
Thus $\E[|N_0 (G_{\tphi_\la}(\Po_\la) ) - \tilde{N}_0|] \to0$,
and (\ref{01claim}) follows.
Note that by (\ref{01claim}) and Markov's inequality,
$\Pr[\tilde{N}_0 (\la) \neq N_0(G_{\phi_\la}(\Po_\la))] \to0 $,
so it suffices to prove (\ref{eq1}) for $\tilde{N}_0(\la)$.

Discretizing space into hypercubes of side $1/m$,
applying the Chen--Stein method of
Poisson approximation and taking the large-$m$ limit as in (32) and
(33) of
\cite{MSTLong} (see also \cite{Penbk}, Theorem~6.7),
we
have that
%
\begin{equation}
\sum_{i=0}^{\infty} \biggl\llvert \Pr\bigl[
\tilde{N}_0(\la) =i\bigr] - \frac{e^{ -\E\tilde{N}_0(\la)}
(\E\tilde{N}_0(\la))^i}{i!} \biggr\rrvert
\leq6(b_1 + b_2), \label{0920a}
\end{equation}
with 
\[
b_1: = \la^2 \int_{\Gamma'_\la} \int
_{B(x;3r_\la^{1-\eps}) \cap
{\Gamma'_\la}} \exp \biggl( - \la\int_{\Gamma} \bigl(
\tilde{\phi}_\la(z-x) + \tilde{\phi}_\la(z-y)\bigr) \,dz
\biggr) \,dy \,dx
\]
and
\begin{eqnarray*}
b_2&:=& \la^2 \int_{\Gamma'_\la} \int
_{B(x; 3 r_\la^{1-\eps}) \cap
{\Gamma'_\la}} \exp \biggl( - \la\int_{\Gamma}
g_{\tilde{\phi}_\la}\bigl(z,\{x,y\}\bigr) \,dz \biggr) \,dy \,dx
\\
&=& 2 \la^2 \int_{\Gamma'_\la} \int_{B(x; 3 r_\la^{1-\eps}) \cap
\Gamma'_{\la,x}}
\exp \biggl( - \la\int_{\Gamma} g_{\tilde{\phi}_\la}\bigl(z,\{x,y\}
\bigr) \,dz \biggr) \,dy \,dx,
\end{eqnarray*}
where for $x \in\Gamma$, if $d=2$, we let $\Gamma'_{\la,x}$ denote
the set of $y \in\Gamma'_{\la}$ lying further from
the boundary of $\Gamma$ than $x $ does, while
if $d \geq3$,
we let $\Gamma'_{\la,x}$ denote
the set of $y \in\Gamma'_\la$ lying closer to
the center of $\Gamma$ in the $\ell_1$ norm than $x$.

By the union bound, $g_{\tilde{\phi}_\la}(z,\{x,y\}) \leq
\tilde{\phi}_\la(z-x)
+ \tilde{\phi}_\la(z-y)$,\vspace*{1pt} and therefore $b_1 \leq b_2$.
Hence by (\ref{0920a})
and (\ref{01claim}),
to prove (\ref{eq1}) it suffices to prove that $b_2 \to0$.

We write $b_2 = b_2^{(1)} + b_2^{(2)}$,
where
$b_2^{(1)}$ denotes the contribution to $b_2$ from
integrating over $(x,y)$ with
$y \in B(x;r_\la)$,
and
$b_2^{(2)}$ denotes the contribution to
$b_2$ from integrating over $(x,y)$ with
$y \in B(x;3 r_{\la}^{1-\eps} ) \setminus B(x;r_\la)$.

First suppose $d=2$.
Using
Lemma~\ref{halflem},
we have that
\[
b_2^{(2)} \leq9 \pi\la^2 r_\la^{2(1-\eps)}
\int_{\Gamma'_\la} \exp \biggl( \biggl( - \la \int
_\Gamma\tilde{\phi}_\la(z-x) \,dz \biggr) -\la(
\eta/4) p_\la r_\la^2 \biggr) \,dx.
\]
By
(\ref{trunc}), we have
%
\begin{equation}
\exp\bigl(- \la I(\phi_\la)\bigr) \leq \exp\bigl(- \la I(\tilde{
\phi}_\la)\bigr) =O\bigl(\la^{-1}\bigr). \label{0725a}
\end{equation}
By (\ref{0729b}),
we have $\exp(- \la p_\la r_\la^2) \leq\exp(-\la
I(\phi_\la)/K(\eta))$, which
is $O( \la^{-1/ K(\eta)} ) $ by (\ref{0725a}). Therefore,
using also (\ref{trunc}) and (\ref{fromThetacond}),
followed by (\ref{0729a}), yields
\[
b_2^{(2)} = O \bigl( \la^{1-\eta/(4K(\eta))}
r_\la^{2(1-\eps)} \bigr) = O \bigl( \la^{\eps- \eta/(4 K(\eta)) } (\log
\la)^{3(1-\eps)} \bigr) \to0.
\]

Now consider $b_2^{(1)}$.
Recall from (\ref{fromThetacond}) that $\la p_\la r_\la^2 = \Theta
(\log
\la)$.
By
Lemma~\ref{halflem},
then (\ref{trunc}) and then (\ref{fromThetacond}),
\begin{eqnarray*}
b_2^{(1)} &\leq& 2 \la^2 \int_{\Gamma'_\la}
\int_0^{r_\la} \exp \biggl( \biggl( - \la\int
_{\Gamma}\tilde{\phi}_\la(z-x) \,dz \biggr) - \la
p_\la(\eta/4) r_\la t \biggr) 2\pi t \,dt \,dx
\\
&=& O \biggl( \la^2 \biggl(\frac{1}{\la} \biggr) \int
_0^\infty \exp\bigl(- (\eta/4) u \bigr) (\la
p_\la r_\la)^{-2} u \,du \biggr) = O \biggl(
\frac{1}{p_\la\log\la} \biggr).
\nonumber
\end{eqnarray*}
Therefore, if $p_\la>1/2$, then $b_2^{(1)} \to0$.
Conversely, if $p_\la\leq1/2$, then
since $g_{\tilde{\phi}_\la}(z,\{x,y\}) \geq\tilde{\phi}_\la(z-x)
+ (1-p_\la)
\tilde{\phi}_\la(z-y)$, and $\phi_\la\in\Phi_{d,\eta}$,
we have
\begin{eqnarray*}
b_{2}^{(1)} &\leq&2 \la^2 \int
_{\Gamma'_\la} \bigl(\pi r_\la^2\bigr) \exp
\biggl( \biggl(- \la \int_{\Gamma} \tilde{\phi}_\la(z-x)
\,dz \biggr) - \la(1-p_\la) \eta p_\la\bigl(\pi
r_\la^2/2\bigr) \biggr) \,dx
\\
&=& O \bigl( \la r_\la^2 \exp\bigl( - \pi(\eta/4) \la
p_\la r_\la^2 \bigr) \bigr)
\end{eqnarray*}
so that by (\ref{0729a}), (\ref{0729b}) and (\ref{0725a})
we have
$
b_{2}^{(1)}
= O  (  ( \log\la )^3
\la^{- \pi\eta/(4K(\eta)) }
 ) = o(1)$.
Hence $b_2^{(1)} \to0$, so that $b_2 \to0$ as required when $d=2$.

Now suppose $d \geq3$.
Let $\tilde{\Gamma}:= \{x \in\Gamma\dvtx  \|x\|_\infty\leq1/2$\}.
Then by Lemma~\ref{halflem}(b),
\begin{eqnarray*}
b_2^{(1)} &\leq&2^{d+1} \lambda^2 \int
_{\tilde{\Gamma}} \int_{B(x;
r_\lambda) \cap
\Gamma'_{\lambda,x}} \exp \biggl( - \lambda
\biggl[\int_\Gamma\tphi_\la(z- x) \,dz\\
&&\hspace*{157pt}{} +
\eta_1 p_\la r_\la^{d-1} |y-x|
\biggr] \biggr) \,dy \,dx
\\
&\leq&2^{d+1} \lambda^2 \int_{\tilde{\Gamma}} \exp
\biggl(- \la\int_\Gamma\tphi_\la(z- x) \,dz \biggr)
\int_{\R^d} \exp\bigl(- \eta_1 \la
p_\la r_\la^d |w| \bigr) r_\la^d
\,dw \,dx,
\end{eqnarray*}
and hence using (\ref{trunc}) followed by (\ref{fromThetacond}),
we obtain that
\[
b^{(1)}_2 = O\bigl( \lambda r_\lambda^d
\bigl( \lambda p_\lambda r_\lambda ^{d}
\bigr)^{-d} \bigr) 
= O\bigl(p_\lambda^{-1}
(\log\lambda)^{1-d}\bigr),
\]
which tends to zero by our assumption on $p_\la$.
By Lemma~\ref{halflem}(b) again,
\[
b_2^{(2)} \leq2^{d+1} \lambda^2 \vold
r_\la^{d(1-\eps)} \int_{\tilde{\Gamma}} \exp \biggl(- \la
\int_\Gamma\tphi_\la(z- x) \,dz \biggr) \times \exp
\bigl( - \eta_1 \la p_\la r_\la^d
\bigr) \,dx,
\]
and hence using (\ref{trunc}),
(\ref{0729a})
and (\ref{fromThetacond}),
with $\delta$ as given at the start of this proof,
we obtain that
$
b^{(2)}_2 = O (
\lambda^{\eps} (\log\la)^{3(1-\eps)}
\exp(- \eta_1 \delta\log\la)
 ).
$
By our choice of $\eps$,
this shows that $b^{(2)}_2 $
tends to zero,
completing the proof.
\end{pf}

\section{Connectivity: The case of small \texorpdfstring{$p_\lambda$}{plambda}}
\label{Sec:psmall}
For any graph $G$,
let $L_2(G)$ denote the order of
its second-largest component, that is, the
second largest of the orders of its components:
if $G$ is connected, set $L_2(G)=0$.
Given the connection functions $(\phi_\la)_{\la>0}$,
let $p_\la$ and $r_\la$ be given by (\ref{plrldef}).
In this section we prove the following result:

\begin{prop}
\label{Poconnthm}
Suppose $(\la(n))_{n \in\N}$ is
an increasing $(0,\infty)$-valued sequence
that tends to $\infty$ as $n \to\infty$,
and for some $\eta\in(0,1]$ and $\alpha\in(0,\infty)$,
$(\phi_\la)_{\la>0}$ is a collection of
connection functions in $\Phi_{d,\eta}$ such that as
$\la\to\infty$ along the sequence $(\la(n))$ we have
(\ref{prcond0}).
Assume
for some $\eps>0$ that
$p_\la= O(\la^{-\eps})$.
Then as $\la\to\infty$ along the same sequence,
\[
\Pr \bigl[ L_2 \bigl(G_{\phi_\la}(\Po_\la)\bigr) >
1 \bigr] \to0.
\]
\end{prop}

It is immediate from Theorem~\ref{Poissecthm}
and Proposition~\ref{Poconnthm}
that under the hypotheses of Proposition~\ref{Poconnthm},
we have
a Poissonized version of (\ref{eq1c}), namely
$
P[ G_{\phi_\la}(\Po_\la) \in\K] \to e^{- \alpha}$.
Our strategy of proof of
Proposition~\ref{Poconnthm} is as follows.
First we shall rule out ``small components'' of
order between 2 and $n^{\eps/2}$ using the Mecke
formula. Then we shall rule out the
possibility of more than one ``large component''
by a ``sprinkling'' argument. That is, we add the edges
in two stages, and even though we make the
number of edges added in the second stage
rather small, with high probability there are enough of them to connect
together any two distinct large components arising
from the first stage.

Given $n \in\N$ and $p \in[0,1]$,
let $G(n,p)$ denote the {Erd\H{o}s--R\'enyi }random graph on
$n$ vertices, that is, the random subgraph of
the complete graph on $n$ vertices, obtained
by including each possible edge independently with
probability $p$.
Our proof of
Proposition~\ref{Poconnthm} uses a lemma on large deviations for
the giant component of
$G(n,p)$.

\begin{lemm}
\label{ERlem}
Suppose $p = p(n)$ is such that
$np \to\infty$ as $n \to\infty$.
Let $E_n$ be the event that
$G(n,p)$ has no component of
order greater than $3n/4$.
Then $\limsup_{n \to\infty} n^{-1} \log P[E_n] <0$.
\end{lemm}

\begin{pf}
Suppose $E_n$
occurs.
Then by starting with the empty set and
adding components of $G(n,p)$
in arbitrary order until we have at least $n/8$ vertices,
we can find a set of between $n/8$ and $7n/8$ vertices that
is disconnected from the rest of the vertices of $G(n,p)$.
Hence by the union bound
and the fact that $e^k \geq k^k/k!$ for any $k$,
\begin{eqnarray*}
P[E_n] &\leq& \sum_{n/8 \leq k \leq7n/8} \pmatrix{n
\cr
k} (1-p)^{k(n-k)} \leq \sum_{n/8 \leq k \leq7n/8}
\frac{n^ke^k}{ k^k} \exp\bigl(-p (7/64)n^2 \bigr)
\\
&\leq& n (8e)^n \exp\bigl(-n^2p /10 \bigr),
\end{eqnarray*}
and the result follows.
\end{pf}

For any graph $G$ any $k \in\N$,
let $T_k(G)$ denote
the number of components of $G$ of order $k$.

\begin{lemm}
\label{lem0711}
Under the hypotheses of Proposition~\ref{Poconnthm},
%
\begin{equation}
\Pr \biggl[ \bigcup_{2 \leq k \leq\la^{\eps/3} } \bigl\{ T_k
\bigl(G_{\phi_\la}(\Po_\la)\bigr) >0 \bigr\} \biggr] \to0.
\label{0711d}
\end{equation}
\end{lemm}

\begin{pf}
We may assume
$r_\la\leq\sqrt{d}$.
By the Mecke formula (\ref{Mecke}), Cayley's formula
(which says there are $k^{k-2}$ trees on $k$ vertices)
and the union bound, $\E T_k( G_{\phi_\la}(\Po_\la))$
is bounded by
\[
\frac{\la^k}{k!} k^{k-2} p_\la^{k-1} \int\cdots
\int \exp \biggl( - \la\int g_{\phi_\la} \bigl(x;\{x_1,
\ldots,x_k\}\bigr) \,dx \biggr) \,dx_1 \cdots \,dx_k
,
\]
where all integrals are over $\Gamma$ in this proof.
By (\ref{0705a}),
this
is bounded by
%
\begin{eqnarray} \label{0621a}
&&\frac{(e \la p_\la)^k}{k^2 p_\la} \int\cdots\int \,dx_1 \cdots \,dx_k
\nonumber
\\[-8pt]
\\[-8pt]
\nonumber
&&\qquad{}\times \exp
\Biggl( - \lambda\int \sum_{i=1}^k
\phi_{\la}(x-x_i) \,dx \Biggr) \exp\bigl(\la k^2
p_\la I(\phi_\la)\bigr).
\end{eqnarray}
By (\ref{Thetacond})
the exponent in the last factor of
(\ref{0621a}) is $O( k^2 p_\la\log\la)$.
If $k \leq\la^{\eps/3} $,
this exponent is $O(1)$, so the
last factor in (\ref{0621a}) is $O(1)$,
uniformly over such $k$.
Thus
\[
\E\sum_{2 \leq k \leq\la^{\eps/3} } T_k \bigl(G_{\phi_\la}(
\Po_\la)\bigr) = O \Biggl( p_\la^{-1} \sum
_{k=2}^\infty\bigl( e p_\la\E
N_0\bigl(G_{\phi_\la}(\Po_\la)\bigr)
\bigr)^k \Biggr),
\]
which tends to zero. Then (\ref{0711d}) follows by
Markov's inequality.
\end{pf}

\begin{pf*}{Proof of Proposition~\ref{Poconnthm}}
Assume
that
$r_\la\leq\sqrt{d}$.
Set $\phi'_\la(x) = \break \phi_\la(x)
(1- \la^{-\eps/6 } ) $ for $x \in\R^d$.
Note that (\ref{prcond0}) still holds
using $\phi'_\la$ instead of $\phi_\la$,
since changing $\phi_\la$ to $\phi'_\la$ gives
an extra term in the exponent of
$O(\la^{1-\eps/6} I(\phi_\la) )$,
which tends to zero by (\ref{Thetacond}).
Also, $\phi'_\la\in\Phi_{d,\eta}$.

Consider generating $G_{\phi_\la}(\Po_\la)$
in two stages. In the first stage, generate $G_{\phi'_\la}(\Po_\la)$.
In the second stage,
for each pair of vertices $X,Y$
not already connected by an edge in the first stage,
add an edge between them with probability
$(\phi_\la(Y-X) - \phi'_\la(Y-X))/(1- \phi'_\la(Y-X))$.

By (\ref{fromThetacond}),
$\la r_\la^d = \Omega(\la^\eps)$ and $r_\la= \Omega(\la^{(\eps-1)/d})$.
We now show that after the first stage,
there is a giant component with high probability.
Partition $\Gamma$ into
cubes of side $1/\lfloor8d/ r_\la\rfloor$.
The number of cubes in the partition
is $O(r_\la^{-d}) = O(\la)$.

By a Chernoff bound (e.g., Lemma~1.2 of \cite{Penbk}),
with high probability
each cube in the partition contains at least
$(9d)^{-d} \la r_\la^d$ vertices
of $\Po_\la$.
Since we assume $r_\la\leq\sqrt{d}$, it is easily verified
that $1/\lfloor8d/r_\la\rfloor\leq r_\la/7d$.
By (\ref{fromThetacond}),
for each cube in the partition,
the restriction of $G_{\phi_\la}(\Po_\la)$
to the vertices within that
cube dominates the {Erd\H{o}s--R\'enyi }random graph $G(n,p)$ with
$np = \Omega( \la r_\la^d (\log\la)/(\la r_\la^d)) =
\Omega( \log\la)$,
so by Lemma~\ref{ERlem},
there is a giant component containing
a proportion of at least $(3/4)$ of the vertices in that cube,
except on an event of
probability $\exp(-\Omega(\la r_\la^d)) =
\exp(-\Omega(\la^{ \eps}))$.
Hence by the union bound,
with high probability the restricted graph within
each of these cubes
contains a giant component.

Also by the same argument,
with high probability, it is the case that for each pair of neighboring
cubes in the partition, the restriction of
$G_{\phi_\la}(\Po_\la)$
to vertices in that pair of cubes
has a giant component with
a proportion of at least $3/4$ of the vertices in
that pair of cubes, and therefore the
two giant components within these neighboring cubes
are connected together. Note that
for any $\delta>0$, with high probability,
by the Chernoff bound,
for each pair of cubes the ratio
of the number of vertices in one cube and the number
of vertices in the other lies between $1-\delta$ and
$1+\delta$.

Hence, after
the first stage there is w.h.p. a giant component
containing a proportion at least $3/4$ of all the
vertices in each of the cubes in the partition.
By Lemma~\ref{lem0711}, also w.h.p. there
is no component of order greater than 1 but less than
$\la^{\eps/3}$.
There may also be some isolated vertices and
some medium-size components of
order between $\la^{\eps/3}$ and $\la/2$.
Now we rule out existence of components
of order greater than $\la^{\eps/3}$
besides the giant component, after the second stage.

After the first stage,
w.h.p.
the giant component contains more than $\lceil
(9d)^{-d} \la r_\la^d/2 \rceil$
vertices in each of the cubes in the partition.
Therefore each vertex not in the
giant component has at least $(9d)^{-d} \la r_\la^d/2$
vertices from the giant component
within the distance of $r_\la$ (viz., those
which are in the same cube of the partition as itself).

Now for each medium-sized component
from the first stage, the probability that it
fails to get attached to the giant component
in the second stage is bounded by
\begin{eqnarray*}
\bigl(1 - \la^{-\eps/6 } \eta p_\la/2 \bigr)^{\la^{\eps/3} \times(9d)^{-d}
\la r_\la^d/2} &\leq&
\exp\bigl( - (9d)^{-d} \eta\la^{ \eps/6 } \la r_\la^d
p_\la/4\bigr)
\\
&\leq& \exp\bigl( - \la^{\eps/6 } \bigr),
\end{eqnarray*}
where the last inequality holds for all large enough $\la$, by
(\ref{fromThetacond}).
The number of medium-sized components from the first stage
is bounded
by $2 \la$ w.h.p., so by the union bound,
the probability that one or more of them
fails to get attached to the giant component tends to zero.

Also the number of isolated vertices from
the first stage is asymptotically Poisson by Lemma~\ref{poislem2},
and the probability that any two of these get connected
together
in the second stage is $O(\la^{-\eps/6} p_\la)$ and thus
tends to zero. Hence w.h.p., after the second
stage there is no component of order greater than 1,
besides the giant component.
\end{pf*}

\section{Connectivity: The case of large \texorpdfstring{$p_\lambda$}{plambda}}
\label{Sec:plarge}
In this section we prove the following result, which extends
Proposition~\ref{Poconnthm} by relaxing the restriction
on $p_\la$ that was imposed there, subject to
$\phi_\la\in\Phi_{d,\eta}^0$.

\begin{theo}
\label{prop0711c}
Let $\alpha\in(0,\infty)$. Suppose that
for some increasing sequence $(\la(n))_{n \in\N}$
that tends to $\infty$ as $n \to\infty$,
$(\phi_\la)_{\la>0}$ satisfies (\ref{prcond0})
as $\la\to\infty$ along the sequence $(\la(n))_{n \in\N}$,
and that there exists $\eta\in(0,1]$ such that
$\phi_\la\in\Phi_{d,\eta}^0$ for all $\la$.
Then
as $\la\to\infty$ along the sequence $(\la(n))_{n \in\N}$,
%
\begin{equation}
\Pr \bigl[ L_2 \bigl(G_{\phi_\la}(\Po_\la)\bigr) >
1 \bigr] \to0. \label{0711f}
\end{equation}
\end{theo}

Throughout this section,
we arbitrarily fix $\eta\in(0,1]$ and assume $\phi_\la\in\Phi
^0_{d,\eta}$
for all $\la>0$, and
$(\phi_\la)_{\la>0}$ satisfy (\ref{prcond0}) for some
$\alpha\in(0,\infty)$ [all asymptotics being as
$\la\to\infty$ along the sequence $(\la(n))_{n \in\N}$].
Define $p_\la:= \mu(\phi_\la)$ and $r_\la:= \rho_\eta(\phi_\la
)$ as
in (\ref{plrldef}), and assume $r_\la=O(1)$.

In view of Proposition~\ref{Poconnthm}, it suffices
to prove the result in
the case where $p_\la= \Omega(\la^{-\eps})$ for some
suitably chosen $\eps>0$. Since the argument is long,
we split the section further by first showing there are no
``small'' components (other than isolated vertices)
and then showing there is not more than one ``large'' component.

\subsection{Small components}

This subsection contains several lemmas
because we sometimes need to distinguish the case with $d=2$
(where we do not assume $\phi_\la\in\Psi_2$), and
we also sometimes distinguish the
case with $p = O(1) $ from $p =o(1)$. Moreover, we distinguish
``very small'' components of (spatial) diameter at most $\delta r_\la$
and ``moderately small'' components of diameter between
$\delta r_\la$ and $(1/\delta) r_\la$, where $\delta$ is
a small (but fixed) constant.

To deal with ``very small'' components (in Lemmas \ref{lem2001a}, \ref
{lem0127},
\ref{lem0624}
and
\ref{lemd3a}) we use the Mecke formula directly and sum over all
possible cardinalities of the component. To deal with
``moderately small components'' (in Lemmas \ref{lemrho1rho2}, \ref{cornerlem}
and \ref{lemd3b}),
we discretize space into
cubes (or strips) of side $\eps r_\la$ for suitably small fixed $\eps$.
For $x \in\Gamma$
and for each possible ``moderately small'' discretized region
(i.e., union of some of
these cubes) containing $x$, we estimate the probability
that the component of
$G_{\phi_\la}(\Po_\la\cup\{x\} )$ containing $x$ is moderately\vadjust{\goodbreak} small
and corresponds to that particular region.
To do this
we show that there is
enough ``unexplored space'' outside the region but inside $\Gamma$,
for the probability of there being no Poisson points in the unexplored
space connected to the cluster within the explored region,
is small compared to the probability of $x$ being isolated.

We need some preliminaries.
First we give
a similar lemma to
Lemma~6 of \cite{MSTLong}.
As before,
let $\bH$ denote the half-space $[0,\infty) \times\R^{d-1}$,
and let $\bQ$ denote the orthant $[0,\infty)^d$.
For $\la>0$
let $\cH_\la^{\bH}:= \cH_\la\cap\bH$, and
let $\cH_\la^{\bQ}:= \cH_\la\cap\bQ$.
Define
\[
\psi_\la(x):= \phi_\la(r_\la x),\qquad x \in
\R^d. 
\]

For any locally finite set $\X$ in $\R^d$, and
any $x \in\R^d$, and connection function $\phi$,
let $C_{\phi}(x,\X) $ be the vertex set of the component
of $G_{\phi}(\X\cup\{x\})$ containing $x$.
Let $D_{\phi}(x,\X):= \diam(C_{\phi}(x,\X)):= \sup_{y,z \in
C_{\phi}(x,\X)}|y-z|$.
For $\A$ a countable set in $\R^2$ and
$x \in\A$,
let $L_{\phi}(x,\A)$ denote the event that
$x$ is the left-most vertex of $C_{\phi}(x,\A)$ (i.e., the first vertex
in the lexicographic ordering).
Also, let $L'_{\phi}(x,\A)$ denote the event that $x$ is the
vertex of $C_{\phi}(x,\A)$ lying closest to the boundary
of the quadrant $\bQ$.

\begin{lemm}
\label{lem2001a}
Suppose $d =2$
and $ p_\la\geq1/2$ for all $\la$.
Then for $0 < \delta\leq\eta/(8\pi)$ we have
\[
\lim_{\la\to\infty} \sup_{x \in\bH} \frac{
P[ 0 < D_{\psi_\la}(x,\cH_{\la r_\la^2}^{\bH}) < \delta; L_{\psi
_\la}(x,
\cH_{\la r_\la^2}^{\bH})]
}{
P[D_{\psi_\la}(x,\cH_{\la r_\la^2}^{\bH}) = 0 ] }
=0.
\]
\end{lemm}

\begin{pf}
Given $x \in\bH$ and $\delta>0$,
let $A_\delta$ denote the right half of the disk of radius $\delta$
centered at $x$.
Let $q_k^\delta(x,\lambda)$ be the probability
that $C_{\psi_\la}(x,\cH_{\la r_\la^2}^{\bH}) $ has precisely $k$ elements
and is contained in $A_\delta$.
Clearly
\[
P\bigl[ 0 < D_{\psi_\la}\bigl(x,\cH_{\la r_\la^2}^{\bH}\bigr) <
\delta; L_{\psi
_\la}\bigl(x,\cH_{\la r_\la^2}^{\bH}\bigr)\bigr]
\leq\sum_{k=2}^\infty q_k^\delta(x,
\la).
\]
By the Mecke formula,
similarly to \cite{conperc}, Proposition~1,
with $h_\phi$ and $g_\phi$ defined at~(\ref{hphidef})
and (\ref{gphidef}),
we have
%
\begin{eqnarray}\label{0920c}
&&q_k^\delta(x,\la)\nonumber \\
&&\qquad= \frac{(\la r_\la^2)^{k-1}}{(k-1)!}
\nonumber
\\[-8pt]
\\[-8pt]
\nonumber
&&\qquad\quad{}\times\int
_{A_\delta} \cdots\int_{A_\delta} h_{\psi_\la}
\bigl(\{ x,x_1,\ldots,x_{k-1}\}\bigr)
\\
&&\hspace*{50pt}\qquad\quad{}\times\exp \biggl( - \la r_\la^2 \int
_{\bH} g_{\psi_\la}\bigl(y,\{x,x_1,
\ldots,x_{k-1}\}\bigr) \,dy \biggr) \,dx_1 \cdots
\,dx_{k-1}. \nonumber
\end{eqnarray}
Similarly $q_1^\delta(x,\la) = \exp(- \la r_\la^2 \int_{\bH} \psi
_\la(y-x) \,dy)$.
Since $h_{\psi_\la}(\A) \leq1$ for any $\A$
we have
%
\begin{eqnarray} \label{0724}
&&\frac{q_k^\delta(x,\la) }{q_1^\delta(x,\la)}\nonumber\\
&&\qquad \leq
\frac{(\la r_\la^2)^{k-1}}{(k-1)!}
\nonumber
\\[-8pt]
\\[-8pt]
\nonumber
&&\qquad\quad{}\times \int_{A_\delta} \cdots
\int_{A_\delta}
 \exp \biggl( - \la r_\la^2 \int
_{\bH} \bigl[g_{\psi_\la}\bigl(y,\{x,x_1,
\ldots,x_{k-1}\} \bigr)\\
&&\qquad\hspace*{185pt}{} - \psi_\la(y-x)\bigr] \,dy \biggr)
\,dx_1\cdots \,dx_{k-1} .\nonumber
\end{eqnarray}
If we restrict the integral in (\ref{0724})
to those $(x_1,\ldots,x_{k-1})$ with $|x_i -x| \leq|x_1 -x|$
for $2 \leq i \leq k-1$,
we reduce it by a factor of $k-1$. Therefore
\begin{eqnarray*}
&&\frac{q_k^\delta(x,\la) }{q_1^\delta(x,\la)}\\
&&\qquad \leq \frac{\la r_\la^2 (\la r_\la^2 \pi/2)^{k-2}}{(k-2)!}\\
&&\qquad\quad{}\times \int_{A_\delta}
|x_1-x|^{2(k-2)}
\\
&&\hspace*{28pt}\qquad\quad{}\times\exp \biggl( - \la r_\la^2 \int
_{\bH} \bigl[g_{\psi_\la}\bigl(y,\{x,x_1\}
\bigr) - \psi_\la(y-x)\bigr] \,dy \biggr)\,dx_1.
\end{eqnarray*}
%
By
Lemma~\ref{halflem}
and the fact that $\rho_\eta(\psi_\la) =1$,
for $x_1 \in A_1$ we have
\[
\int_{\bH} \bigl[g_{\psi_\la}\bigl(y,\{x,x_1
\}\bigr) - \psi_\la(y-x)\bigr] \,dy 
\geq |x_1-x| \eta p_\la/4,
\]
so that for $\delta\leq1$ we have
\[
\frac{q_k^\delta(x,\la) }{q_1^\delta(x,\la)} \leq \frac{\la r_\la^2 (\la r_\la^2 \pi/2)^{k-2}}{(k-2)!} \int_{A_\delta}
|x_1-x|^{2(k-2)} \exp \bigl( - \la r_\la^2
(\eta/4) p_\la |x_1-x| \bigr) \,dx_1 .
\]
Summing over $k \geq2$ and using the assumptions
$p_\la\geq1/2$ and $\delta\leq\eta/ (8 \pi)$, yields
\begin{eqnarray*}
\sum_{k=2}^\infty \frac{q_k^\delta(x,\la) }{q_1^\delta(x,\la)} &\leq&
\la r_\la^2 \int_{A_\delta} \exp\bigl( \la
r_\la^2 \bigl[(\pi/2) |x_1-x|^2 -
(\eta/4) p_\la|x_1 -x|\bigr] \bigr) \,dx_1
\\
&\leq&\la r_\la^2 \int_{A_\delta} \exp
\bigl( - \la r_\la^2 |x_1 -x| \eta/16\bigr)
\,dx_1 = O\bigl(\bigl(\la r_\la^2
\bigr)^{-1}\bigr),
\end{eqnarray*}
which tends to zero by (\ref{fromThetacond}).
\end{pf}

In the case with $p_\la\leq1/2 $,
we give a similar result to the last one, but
for general $d \geq2$.
Let $\vold$ denote
the volume of the unit ball in $d$ dimensions.
Let $\tilde{\bQ}$ denote the orthant
$\bQ$ if $d \geq3$, but denote the half-space $\bH$
if $d=2$.

\begin{lemm}
\label{lem0127}
Suppose $\phi_\la$ and $\psi_\la$ are as before (now for general $d$,
$d \geq2$).
Let $ 0 < \delta< \eta/8$.
If $p_\la\leq1/2$ for all $\la$ but
$p_\la= \Omega( \la^{- 1/2^{d+3}})$,
then
\[
\lim_{\la\to\infty} \sup_{x \in\tilde{\bQ}} \biggl(
\frac{ P[0 < D_{\psi_\la}(x,\cH^{\tilde{\bQ}}_{\la r_\la^d})
< \delta]}{P [ D_{\psi_\la}(x,\cH^{\tilde{\bQ}}_{\la r_\la^d})
=0] } \biggr) =0.
\]
\end{lemm}
\begin{pf}
For $\delta>0$, $x \in\tilde{\bQ}$ and $k \in\N$, define
\[
w_\la(k,\delta) := \frac{ P[0 < D_{\psi_\la}(x,\cH^{\tilde{\bQ}}_{\la r_\la^d})
< \delta;
\card(C_{\psi_\la}( x, \cH^{\tilde{\bQ}}_{\la r_\la^d})) =k+1 ]}{
P [ D_{\psi_\la}(x,\cH^{\tilde{\bQ}}_{\la r_\la^d}) =0] },
\]
where $\card(\cdot)$ denotes the number of elements in a set.
For $k \in\N$ we have,
similarly to
(\ref{0920c}),
that
%
\begin{eqnarray}\label{0127b}
 w_\la(k,\delta)&\leq &\frac{(\la r_\la^d)^{k}}{k!}\nonumber \\
&&{}\times\int_{B(x;\delta) \cap\tilde{\bQ}}
\cdots\int_{B(x;\delta)
\cap\tilde{\bQ}}
\nonumber
\\[-8pt]
\\[-8pt]
\nonumber
&&{}\times\exp \biggl( - \int_{\tilde{\bQ}} \la r_\la^d
\bigl[g_{\psi_\la}\bigl(y,\{x,x_1,\ldots,x_{k}\}
\bigr)\\
&&\hspace*{118pt}{}- \psi_\la(y-x)\bigr] \,dy \biggr) \,dx_1 \cdots
\,dx_{k} .\nonumber
\end{eqnarray}
Now,
%
\begin{eqnarray}\label{0723c}
&&g_{\psi_{\la}}\bigl(y,\{x,x_1,x_2,
\ldots,x_k\}\bigr) -\psi_\la( y-x)\nonumber \\
&&\qquad\geq
(1-p_\la) \Biggl( 1 - \prod_{i=1}^k
\bigl( 1- \psi_\la(y-x_i)\bigr) \Biggr)
\\
&&\qquad\geq (1-p_\la) \Biggl(1 - \exp \Biggl(- \sum
_{i=1}^k \psi_\la (y-x_i)
\Biggr) \Biggr). \nonumber
\end{eqnarray}
First consider $k \leq1/p_\la$.
Since $1 -e^{-x} \geq x/2$ for $0 \leq x \leq1$,
and we assume $p_\la\leq1/2$,
for such $k$
we have
\[
g_{\psi_{\la}}\bigl(y,\{x,x_1,x_2,
\ldots,x_k\}\bigr) -\psi_\la( y-x) \geq(1/4) \sum
_{i=1}^k \psi_\la(y-x_i).
\]
Now $\int_{\tilde{\bQ}} \psi_\la(y-x_i) \,dy \geq I(\psi_\la)/2^d$
for each $\la$ and each $x_i$
because
for $d \geq3$
we assume $\phi_\la\in\Psi_d$, and
for $d=2$ we assume
$\tilde{\bQ} = \bH$,
and
$\phi_\la$
satisfies $\phi_\la(x) =\phi_\la(-x)$ for all $x $.
Therefore by (\ref{0127b}),
for $k \leq1/p_\la$ we have
\begin{eqnarray*}
w_\la(k,\delta) &\leq& \frac{(\la r_\la^d)^{k}}{k!} 
\int
_{(B(x;\delta) \cap\tilde{\bQ})^d} \exp \Biggl( - \frac{1}{4} \int
_{\tilde{\bQ}} \la r_\la^d \sum
_{i=1}^k \psi_\la(y-x_i) \,dy
\Biggr) 
\,d(x_1, \ldots, x_{k} )
\\
&\leq&\frac{(\delta^d \vold\la r_\la^d)^{k}}{k!} \exp\bigl( - \la k I(\phi_\la)/2^{d+2}
\bigr).
\end{eqnarray*}
Hence
%
\begin{equation}
\sum_{k=1}^{\lfloor1/p_\la\rfloor} w_\la(k,
\delta) \leq \exp\bigl[ \delta^d \pi_d \la
r_\la^d e^{- \la I(\phi_\la)/2^{d+2} }\bigr] -1. \label{0127a}
\end{equation}
Since we assume (\ref{prcond0}) we have
$e^{-\la I(\phi_\la)} = O(\la^{-1})$,
and using (\ref{fromThetacond}) we have that
\[
\la r_\la^d e^{- \la I(\phi_\la)/2^{d+2}} = O \biggl(
\frac{ \log\la}{ p_\la\la^{1/2^{d+2}}} \biggr),
\]
which tends to zero,
by our condition on $p_\la$. Therefore the expression in
(\ref{0127a}) tends to zero.

Now consider $k > 1/p_\la$.
For
$x_1,\ldots,x_k \in B(x;\delta)$ and
$y \in B(x;1/2)$ we have
$|y-x_i| < 1$ and hence $\psi_\la(y-x_i) \geq\eta p_\la$
for $1 \leq i \leq k$,
so by (\ref{0723c})
we have that
\begin{eqnarray*}
g_{\psi_{\la}}\bigl(y,\{x,x_1,x_2,
\ldots,x_k\}\bigr) -\psi_\la( y-x) &\geq&(1-p_\la)
\bigl(1- \exp( - \eta k p_\la)\bigr)\\
& \geq&\bigl(1- e^{-\eta}
\bigr)/2.
\end{eqnarray*}
Therefore using (\ref{0127b})
and the fact that $1-e^{-\eta} \geq\eta/2$,
we have
\[
\sum_{k > 1/p_\la} w_\la(k,\delta) \leq \exp
\bigl( \delta^d \vold\la r_\la^d\bigr) \exp
\bigl(-\vold\la r_\la^d \eta/2^{d+2}\bigr),
\]
and by the choice of $\delta$, this
tends to zero. Combining these estimates gives
the result.
\end{pf}

Combining Lemma~\ref{lem2001a} and
the case $d=2$ of Lemma~\ref{lem0127}
immediately gives us the following.

\begin{lemm}
\label{lem0624}
Suppose $d=2$ and
$\eta$, $\phi_\la$ and $\psi_\la$ are as before. Suppose
also that
$p_\la= \Omega( \la^{ -1/32 })$.
Let $0 < \delta< \eta/(8\pi)$.
Then
\[
\lim_{\la\to\infty} \sup_{x \in\bH} \biggl(
\frac{ P[0 < D_{\psi_\la}(x,\cH^{\bH}_{\la r_\la^2})
< \delta; L_{\psi_\la}(x;\cH^{\bH}_{\la r_\la^2})]}{P [ D_{\psi
_\la}(x,\cH^{\bH}_{\la r_\la^2}) =0] } \biggr) =0.
\]
\end{lemm}

\begin{lemm}
\label{lemrho1rho2}
Given $0< \delta< \rho< \infty$, it is the case
(for general $d \geq2$) that
\[
\lim_{\la\to\infty} \sup_{x \in\bH} \frac{P[ \delta
< D_{\psi_\la}(x,\cH_{\la r_\la^2}^{\bH}) < \rho; L_{\psi_\la
}(x,\cH_{\la r_\la^2}^{\bH})]}{
P[D_{\psi_\la}(x,\cH_{{\la r_\la^2}}^{\bH}) = 0 ] }
=0.
\]
\end{lemm}
\begin{pf}
This can be proved along the lines of
\cite{conperc}, Lemma~3;
the argument still works in the case with $p_\la\to0$,
provided $\la r_\la^2 p_\la\to\infty$, which is always the
case by (\ref{fromThetacond}).
\end{pf}

Similarly to
\cite{MSTLong}, Lemma~7 (which is missing a factor
of $\pi$ in the exponent) we have the following:

\begin{lemm}
\label{cornerlem}
Suppose $d=2$.
For any $\rho>0$,
as $\la\to\infty$ we have
%
\begin{equation}
%
\sup
_{x \in\bQ} P\bigl[D_{\psi_\la}\bigl(x; \cH_{\la r_\la^2}^\bQ
\bigr) < \rho\bigr] = o \bigl( \exp\bigl\{ - \la\eta I(\phi_\la) /
\bigl(3 K(\eta)\bigr) \bigr\} \bigr) . \label{0123a}
\end{equation}
\end{lemm}

\begin{pf}
Fix $\rho>0$.
Divide $\bQ$ into vertical strips of width $1/9$, denoted
$S_i, i \in\N$, where $S_i:= [(i-1)/9,i/9) \times[0,\infty)$.
Let $x \in\bQ$, and let $i_0 = i_0(x)$ be
the choice of $i$ such that $x \in S_i$.
Also let $i_1 = i_0 + 9 \lceil\rho\rceil$.

Given $\la$,
for $i \in\N\cap[i_0,i_1]$
let $E'_{i}$ be the event that
the right-most point of
$C_{\psi_\la}(x,\cH_{\la r_\la^2} ^\bQ) $
lies in $S_i$. If
$D_{\psi_\la}(x; \cH_{\la r_\la^2}^\bQ) < \rho$, then
one of the events $E_{i_0}, \ldots, E_{i_1}$ occurs.

Now fix $i \in\N\cap[i_0,i_1]$. Set
$A_i:= \bigcup_{j \leq i} S_j$ and $A_i^c:= \bigcup_{j >i} S_j$.
Consider
generating $G_{\psi_\lambda}( \{x\} \cup\cH_{\la r_\la^2}^\bQ)$
in two stages. In the first stage, generate the Poisson
process $\cH_{\la r_\la^2} \cap A_i$,
and add edges between points of
$ \{x \} \cup(\cH_{\la r_\la^2} \cap A_i)$ with
probabilities determined by
the connection function $\psi_\la$. Then in the
second stage, add the points of $\cH_{\la r_\la^2} \cap A_i^c $,
and add edges between these added
points, and between the added points and the
points from the first stage, again
using the connection function $\psi_\la$.

The first stage generates a realization of the graph
$G_{\psi_\lambda}( \{x\} \cup(\cH_{\la r_\la^2} \cap A_i) )$:
let $E_{i,1}$ be the event that the resulting realization of
$C_{\psi_\la}(x,\cH_{\la r_\la^2} \cap A_i) $
includes at least one vertex in $S_i$.
Let $E_{i,2}$ be the event that the second stage does not
generate any new Poisson points that are connected to
vertices of
$C_{\psi_\la}(x,\cH_{\la r_\la^2} \cap A_i) $ arising
from the first stage. Then $E'_i = E_{i,1} \cap E_{i,2}$.

Suppose $E_{i,1}$ occurs. Let $z$ be the right-most vertex
of $C_{\psi_\la}(x,\cH_{\la r_\la^2} \cap A_i) $;
then $z \in S_i$ by definition.
Then in stage 2, a necessary condition
for $E_{i,2}$ to occur is that
there is no point of
$\cH_{\la r_\la^2} \cap A_i^c $ connected by an edge
$z$. Since $B(z;1) \cap A_i^c$
has area of at least $(\pi/4) - 1/9$,
\[
\Pr\bigl[ E'_i | E_{i,1} \bigr]\leq\exp
\bigl( - \la r_\la^2 \eta p_\la\bigl( (\pi/4)
- 1/9\bigr)\bigr) 
\leq\exp\bigl( - \eta\lambda I(
\phi_\la) /\bigl(2K(\eta)\bigr) \bigr),
\]
where the last inequality comes from
(\ref{0729b}). This gives us (\ref{0123a}).
\end{pf}

For $x \in\Gamma$, let $\Gamma_x$ be the set of $y \in\Gamma$
such that $y$ is closer to the center of
$\Gamma$ in the $\ell_1$ norm than $x$.
For $\rho>0$ and $x \in\Gamma$, let
$E_{\la,\rho,x}$ be the event that there is a nonempty
set $U$ of
points of $\Po_\la$ contained in $B(x;\rho) \cap\Gamma_x$,
such that no other point of $\Po_\la\setminus U$ is connected to any
point of
$\{x \} \cup U$ in $G_{\phi_\la}(\{x\} \cup\Po_\la)$.

\begin{lemm}
\label{lemd3a}
Suppose $d \geq3$
and $p_\la=\Omega(1)$.
Then there exists $\delta>0$ such that
\[
\lim_{\la\to\infty} \sup_{x \in\Gamma} P[
E_{\la, \delta r_\la,x}] \Bigm/ \exp \biggl(- \la\int_{\Gamma} \phi
_\la(y-x) \,dy \biggr) =0.
\]
\end{lemm}

\begin{pf} The proof resembles that
of \cite{kcon}, Lemma~5.2 or \cite{Penbk}, Lemma~13.15.
For $j \in\N$ let
$\mu^x(j,\la)$ be the number of subsets $U$ of
$\Po_\la$ with $j$ elements, such
that $U \subset\Gamma_x \cap B(x;\delta r_\la)$, and no
element of $U \cup\{x\}$ is connected to any
element of $\Po_\la\setminus U$
in $G_{\phi_\la}(\{x\} \cup\Po_\la)$.
Then by the Mecke formula (\ref{Mecke}),
\begin{eqnarray*}
\E\mu^j(x,\la) &= &\frac{\la^j}{(j-1)!} \int_{\Gamma_x \cap B(x;\delta r_\la) }
\int_{(\Gamma_x \cap B(x;|y-x|) )^{j-1} }
\\
&&{}\times \exp \biggl( -\la\int g_{\phi_\la}\bigl(z,\{x,y,x_1,
\ldots,x_{j-1}\}\bigr) \,dz \biggr) \,d(x_1,\ldots,x_{j-1})
\,dy
\\
&\leq&\frac{ \la( \la\pi_d )^{j-1} }{ (j-1)! }
\int_{\Gamma_x \cap B(x;\delta r_\la) } |y-x|^{d(j-1)}\\
&&{}\times \exp
\biggl( - \la\int g_{\phi_\la}\bigl(z,\{x,y\} \bigr) \,dz \biggr) \,dy.
\end{eqnarray*}
Assume $\delta\leq1$.
By Lemma~\ref{halflem}(b), the integrand in
the last exponent is bounded below by $\phi_\la(z-x)
+ \eta_1 p_\la\rho_\la^{d-1} |y-x|$, and therefore
\begin{eqnarray*}
&&\frac{ \E\mu^j(x,\la) }{
\exp( - \la\int\phi_\la(z-x) \,dz ) } \\
&&\qquad\leq \frac{ \la( \la\pi_d )^{j-1} }{ (j-1)! } \int_{B(x;\delta r_\la)}
|y-x|^{d(j-1)} \exp\bigl( - \eta_1 \la p_\la
r_\la^{d-1}|y-x| \bigr)\,dy.
\end{eqnarray*}
Summing over $j $ and changing variable to $w = (y-x)/r_\la$, we
obtain
\begin{eqnarray*}
\frac{ P[ E_{\la, \delta r_\la,x}] }{
\exp (- \la\int_{\Gamma} \phi_\la(z-x) \,dz
 ) } &\leq&\la\int_{B(x;\delta r_\la)} \exp\bigl( \la
\pi_d |y-x|^d - \eta_1 \la p_\la
r_\la^{d-1}|y-x| \bigr) \,dy
\\
&=& \la\int_{B(0;\delta)} \exp\bigl( \la\pi_d
r_\la^d |w|^d - \eta_1 \la
p_\la r_\la^d |w| \bigr) r_\la^d
\,dw.
\end{eqnarray*}
Using our assumption
on $p_\la$,
we may choose $\delta$ small enough so that
$\pi_d \delta^d \leq(\eta_1/2) p_\la\delta$ for all $\la$,
and then there is a
constant $\delta'$\vspace*{1pt} so the last bound is at most
$\la r_\la^d \int\exp(- \delta' \la r_\la^d |w|) \,dw$,
which is $O((\la r_\la^d)^{1-d})$, and therefore tends to zero
by~(\ref{fromThetacond}).
\end{pf}

In the next lemma we do not need to assume $p_\la= \Omega(1)$.

\begin{lemm}
\label{lemd3b}
Suppose $d \geq3$.
Then for $0 < \delta< \rho< \infty$ we have
\[
\lim_{\la\to\infty} \sup_{x \in\Gamma} P[ E_{\la, \rho r_\la,x}
\setminus E_{\la,\delta r_\la,x} ]\Bigm / \exp \biggl(- \la\int_{\Gamma}
\phi_\la(y-x) \,dy \biggr) =0.
\]
\end{lemm}

\begin{pf}
Fix $\delta$ and $\rho$, and assume $\delta\leq1$.
Let $\eps>0$ be a small constant to be chosen later.
Given $\la$,
divide $\R^d$ into boxes [i.e., hypercubes of the form
$\prod_{i=1}^d[a_i,a_i+h)$]
of side $h= \eps r_\la$. Let $\Lambda'_\la$ be the set
of centers of these boxes. For $z \in\Lambda'_\la$ let
$B'_z$ be the box centered at $z$.
Let $x \in\Gamma$, and
let $z_x $ be the $z \in\Lambda'_\la$ such that $x $ lies
in $B'_z$.
Also, for all $z \in\Lambda'_\la$ let $B_z:= B'_z \cap\Gamma_x$.

For $\sigma\subset\Lambda'_\la$, let
$B_\sigma:= \bigcup_{z \in\sigma} B_z $.
Let ${\cal C}(\la,x)$ be the set of $\sigma\subset\Lambda'_\la$
such that: (i)~$z_x \in\sigma$,
(ii)~$\sigma\subset B(x; (\rho+ d \eps) r_\la)$, (iii)~$\sigma\setminus B(x;(\delta- d \eps) r_\la) \neq\varnothing$ and
(iv)~$|B_z| >0 $ for each $z \in\sigma$ (where $|\cdot|$ denotes Lebesgue
measure).
In the sequel, we assume $\eps< \delta/(2d)$ so
that $\delta- d \eps> \delta/2$.

For $\sigma\in{\cal C}(\la,x)$, let $E'_\la(\sigma)$ be the event
that: (i) $ \sigma= \{z \in\Lambda'_\la\dvtx  C_{\phi_\la}(x,\Po_\la)
\cap B_z \neq\varnothing\}$
and (ii) $C_{\phi_\la}(x,\Po_\la) \subset\Gamma_x$.
Then $ E_{\la, \rho r_\la,x} \setminus E_{\la,\delta r_\la,x}
\subset\bigcup_{\sigma\in{\cal C}(\la, x) } E'_\la(\sigma)$.

Let $\sigma\in{\cal C} (\la, x)$.
Consider generating $C_{\phi_\la}(x,\Po_\la)$ in two stages,
similarly to the proof of Lemma~\ref{cornerlem}. In stage 1,
add all points of $\Po_\la$ in $B_\sigma$, all edges
involving these points and $x$ (using the connection function
$\phi_\la$). In stage 2, add the points of $\Po_\la$
in $\Gamma\setminus B_\sigma$,
add connections between these new points and each other
and between the new points and the points from the stage 1,
again using connection function $\phi_\la$.

In stage 1, we generate a realization of $G_{\phi_\la}(\{x\} \cup
(\Po_\la
\cap B_\sigma))$, and hence a
realization of $C_{\phi_\la}(x,\Po_\la\cap B_\sigma)$.
Let $E'_{\la,1}(\sigma)$ be the event that this
realization of $C_{\phi_\la}(x, \Po_\la\cap B_\sigma)$
is contained in $\Gamma_x$ and
includes at least one point from each $B_z,z \in\sigma$.
Let $E'_{\la,2}$ be the event that
none of the new points created in stage 2 are joined
to any points of the realization of
$C_{\phi_\la}(x, \Po_\la\cap B_\sigma)$ generated in stage~1.
Then $E'_\la(\sigma) = E'_{\la,1}(\sigma) \cap E'_{\la,2}(\sigma)$.
Since the cardinality of ${\cal C}(\la,x)$ is
bounded independently of $x$ and $\la$,
it suffices to show that
%
\begin{equation}
\limsup_{\la\to\infty} \sup_{x \in\Gamma, \sigma\in{\cal C}(\la,x) }
\frac{ \Pr[E'_{\la,2}(\sigma)|E'_{\la,1}(\sigma)] }{\exp(- \la
\int
\phi_\la(y-x) \,dy ) } =0. \label{0116a}
\end{equation}
Now,
\[
\Pr\bigl[E'_{\la,2}(\sigma)|E'_{\la,1}(
\sigma)\bigr] \leq\exp \biggl( - \la\inf_{\cX\subset B_\sigma\cap\Gamma_x\dvtx \cX\cap Q_z \neq\varnothing\ \forall z \in\sigma} \int
_{\Gamma\setminus B_\sigma} g_{\phi_\la} (y;\cX) \,dy \biggr),
\]
and for each $\cX\subset B_\sigma\cap\Gamma_x$
with $\cX\cap Q_z \neq\varnothing$ for all $ z \in\sigma$,
we have
\begin{eqnarray*}
\int_{\Gamma\setminus B_\sigma} g_{\phi_\la} (y;\cX) \,dy 
&=&
p_\la\int_0^{p_\la^{-1}} \int
_{\Gamma\setminus B_\sigma} {\mathbf1}_{
\{ g_{\phi_\la}(y; \X) \geq p_\la u \} } \,dy \,du
\\
&\geq& p_\la\int_0^1 \bigl| \Gamma\cap
\bigl(\sigma\oplus B\bigl(0;\rho_u (\phi_\la) - d \eps
r_\la\bigr) \bigr) \setminus B_\sigma\bigr| \,du,
\end{eqnarray*}
where
the last line arises because if $y \in\sigma\oplus B(0;\rho_u( \phi
_\la)
- d \eps r_\la)$, then there exists $v \in\cX$
with $|y-v| \leq\rho_u(\phi_\la)$
and therefore $g_{\phi_\la}(y;\cX) \geq\phi_\la(y-v) \geq u p_\la$
by (\ref{rhodef}).

For $0 < u \leq1$, since $\phi_\la\in\Phi_{d,\eta}^0$ we have
$r_\la\leq\rho_u (\phi_\la) \leq\eta^{-1}r_\la$.
Hence $\rho_u(\phi_\la) - d \eps r_\la\geq r_\la/2$.
Using \cite{kcon}, Proposition~2.1 or
\cite{Penbk}, Proposition~5.15 and writing
$V_r(x)$ for $|B(x;r) \cap\Gamma|$, we can find
a constant $\eta_3 >0$,
depending only on $d$ and~$\eta$, such that
\[
\int_{\Gamma\setminus B_\sigma} g_{\phi_\la} (y;\cX) \,dy \geq
p_\la \biggl(\int_0^1
V_{\rho_u (\phi_\la) - d \eps r_\la}(x) \,du + \int_0^1
\eta_3 r_\la^d \,du \biggr).
\]
%
To estimate the first term in the expression above, note that
since $\rho_u(\phi_\la) \leq\eta^{-1} r_\la$, there is a constant $K_1$
(depending on $d$ and $\eta$) such that
$
V_{\rho_u(\phi_\la)}(x) \,du
- V_{\rho_u (\phi_\la) - d \eps r_\la}(x)
\leq K_1 r_\la^d \eps$.
Therefore
\begin{eqnarray*}
\int_{\Gamma\setminus B_\sigma} g_{\phi_\la} (y;\cX) \,dy &\geq&
p_\la \biggl( \int_0^1
V_{\rho_u(\phi_\la)}(x) \,du - K_1 r_\la^d \eps+
\eta_3 r_\la^d \biggr)
\\
&=& \int_{\Gamma} \phi_\la(y-x) \,dy - K_1
p_\la r_\la^d \eps+ \eta_3
p_\la r_\la^d,
\end{eqnarray*}
and by choosing $\eps< \eta_3/(2K_1) $,
we have
that the ratio on the left-hand side of (\ref{0116a}) is bounded below by
$
\exp(- \eta_3 \lambda p_\la r_\la^d/2)$,
uniformly over $x$ and $\sigma$.
Since $\la r_\la^d p_\la\to\infty$
by (\ref{fromThetacond}),
this gives us (\ref{0116a}) as required.
%
\end{pf}

Given $\la>0$, $\rho>0$, define the event
\[
E_{\la}^{\rho} = \bigl\{ \exists x \in\Po_\la\dvtx
0 < D_{\phi_\la} (x,\Po_\la) \leq\rho\bigr\}.
\]

\begin{prop}
\label{lemd3c}
Let $\eta\in(0,1]$, $\alpha\in(0,\infty)$ and
$0 < \eps\leq\min( \eta/ \break (7K(\eta)),  2^{-(d+3)} )$.
Suppose $\phi_\la\in\Phi_{d,\eta}^0$ for all $\la$,
(\ref{prcond0}) holds, and
$p_\la= \Omega(\la^{-\eps})$.
Then for any $\rho>0$, we have
$\lim_{\la\to\infty} P[E^{\rho r_\la}_\la] =0$.
\end{prop}

\begin{pf}
First consider the case with $d \geq3$.
Assume first that
$p_\la=\Omega(1)$. Then by the Mecke formula and
the preceding two lemmas,
\[
P \bigl[ E_\la^{\rho r_\la} \bigr] \leq\int_\Gamma
P[ E_{\la,\rho r_\la,x} ]\la \,dx = o(1) \times\int_{\Gamma} \exp
\biggl( - \la\int_{\Gamma} \phi_\la(y-x) \,dy \biggr)
\la \,dx,
\]
which is $o(1)$ by (\ref{prcond0}).

Now suppose instead that $p_\la\to0$ but
$p_\la= \Omega(\la^{-\eps})$.
Then $r_\la= o(1)$ by (\ref{fromThetacond}).
Let $\tilde{\Gamma}$ denote the set of points
in $\Gamma$ lying closer to the origin (in the Euclidean norm) than
to any other corner of $\Gamma$.
Choosing $\delta\in(0, \eta/8)$
we have by the Mecke formula and
Lemma~\ref{lem0127}
that
\begin{eqnarray*}
\Pr\bigl[ E_\la^{\delta r_\la} \bigr] &\leq& 2^d \la \int
_{\tilde{\Gamma}} P\bigl[ 0 < D_{\phi_\la}(x, \Po_\la) <
\delta r_\la\bigr] \,dx
\\
&= &2^d \la \int_{\tilde{\Gamma}} P\bigl[ 0 <
D_{\psi_\la}\bigl(r_\la^{-1}x, \cH_{\la r_\la^d}^{\bQ}
\bigr) < \delta\bigr] \,dx
\\
&=& o(1) \times \la\int_{\tilde{\Gamma}} P\bigl[ D_{\psi_\la}
\bigl(r_\la^{-1}x, \cH_{\la r_\la^d}^{\bQ}\bigr)
=0 \bigr] \,dx,
\end{eqnarray*}
which tends to zero by
(\ref{prcond0}).
Also for any finite $\rho> \delta$,
by the Mecke formula,
\[
\Pr\bigl[ E_{\la}^{\rho r_\la} \setminus E_\la^{\delta r_\la}
\bigr] \leq\la\int_\Gamma \Pr[ E_{\la,\rho r_\la, x} \setminus
E_{\la,\delta r_\la,x} ] \,dx,
\]
which tends to zero
by Lemma~\ref{lemd3b} and
(\ref{prcond0}).
This gives us the result for the case with $d \geq3$.

Now consider the case with $d =2$.
Then
$r_\la^2 = O(\la^{2 \eps-1})$
by (\ref{fromThetacond}).
Let $T_1$ (resp., $T_2$, $T_3,T_4$)
be the set of points of $[0,1]^2$ that
lie closer to the left (resp., top, right, bottom)
edge of $\Gamma$ than to any of the other edges of $\Gamma$
[so $T_1$ is the triangle with corners
at $(0,0),(0,1)$ and $(1/2,1/2)$].

For $x \in\Gamma$, let $\tilde{L}_{\phi_\la}(x,\Po_\la)$
be the event that $x$ is the point of $C_{\phi_\la}(x, \Po_\la)$
lying closest to the
boundary of $[0,1]^2$.
Let $M_\la$ be the number of $x \in\Po_\la$ such that
(i)~$D_{\phi_\la}(x,\Po_\la) < \rho r_\la$,
and (ii)~$x$ is the point of $C_{\phi_\la}(x, \Po_\la)$
nearest to the boundary of $\Gamma$.
Then by the Mecke equation,
\[
P\bigl[E_\la^{\rho r_\la}\bigr] \leq\E M_\la= \sum
_{i=1}^4 a_i,
\]
where we set
\[
a_i:= \la\int_{T_i} P\bigl[ 0 <
D_{\phi_\la}(x, \Po_\la) < \rho r_\la;
\tilde{L}_{\phi_\la}(x,\Po_\la)\bigr] \,dx.
\]
We consider just $a_1$ (the other terms are treated similarly).
Let $T_{1,1}$ be the part of $T_1$ away from the corner
of $\Gamma$, defined by
\[
T_{1,1}:= T_1 \setminus\bigl( \bigl[0,2 \bigl(\rho+
\eta^{-1}\bigr) r_\la\bigr] \times \bigl( \bigl[0,2 \bigl(
\rho+ \eta^{-1}\bigr) r_\la\bigr] \cup\bigl[1-2 \bigl(\rho+
\eta^{-1}\bigr) r_\la,1\bigr]\bigr)\bigr).
\]
Let $a_{1,1}$ be the contribution to $a_1$ from $x \in T_{1,1}$.
Using our assumption that $\phi_\la\in\Phi_{d,\eta}^0$,
we have
\begin{eqnarray*}
a_{1,1} &\leq&\la \int_{T_{1,1}} P\bigl[ 0 <
D_{\phi_\la} (x, \Po_\la) < \rho r_\la;
L_{\phi_\la}(x,\Po_\la) \bigr] \,dx
\\
&=& \la \int_{T_{1,1}} P\bigl[ 0 < D_{\psi_\la}
\bigl(r_\la^{-1}x, \cH_{\la r_\la^2}^{\bH}\bigr)
< \rho; L_{\psi_\la}\bigl(r_\la^{-1}x,
\cH^{\bH}_{\la r_\la^2}\bigr) \bigr] \,dx.
\end{eqnarray*}

Now using Lemmas \ref{lem0624} and
\ref{lemrho1rho2}
we obtain that
\begin{eqnarray*}
a_{1,1} &=& o(1) \times\int_{T_{1,1}} \la\Pr\bigl[
D_{\psi_\la}\bigl(r_\la^{-1} x,\cH_{\la r_\la^2}^{\bH}
\bigr) = 0 \bigr] \,dx\\
& =& o(1) \times\la\int_{T_{1,1}} \Pr\bigl[
D_{\phi_\la}( x,\Po_{\la}) = 0 \bigr] \,dx,
\end{eqnarray*}
which tends to zero by (\ref{prcond0}).

Let $a_{1,2}$ be the contribution to $a_{1}$
from $x \in T_1 \cap[0,2 (\rho+ \eta^{-1}) r_\la]^2$.
By Lemma~\ref{cornerlem},
\[
a_{1,2} \leq\la\bigl( 2 \eta^{-1} r_\la
\bigr)^2 \exp\bigl(- \la\eta I(\phi_\la) /\bigl(3 K(\eta)
\bigr)\bigr) = O\bigl( \la^{2 \eps-\eta/(3K(\eta))}\bigr),
\]
where for the last estimate we use (\ref{fromThetacond})
and (\ref{prcond0}).
Thus $P[E_\la^{\rho r_\la}] \to0$.
\end{pf}

\subsection{Large components}
\label{Subsec:LC}

In this section we implement the strategy mentioned in
the final paragraph of Section~\ref{secmainres}.
In the sequel, given
$\la>0$ we couple the graphs
$G_{\phi_\la}(\Po_\la\cap A), A \subset\R^d$ in the following,
natural way. For $A \subset\R^d$ we define
$G_{\phi_\la}(\Po_\la\cap A)$ to be the subgraph
of $G_{\phi_\la}(\Po_\la)$ induced by the vertex set
$\Po_\la\cap A$.

Given $\la$, let $m_\la: = \lceil2d /r_\la\rceil$.
Set $\Lambda_\la:= \{0,1,\ldots, m_\la-1\}^d$.
For $z \in\Lambda_\la$ let
$Q_z$ denote the cube $\{m_{\la}^{-1}z \} \oplus[0,1/m_\la)^d$,
and let $\overline{Q}_z$ denote the closure of $Q_z$.
The cubes $Q_z, z \in\Lambda_\la$
form a partition of $[0,1)^d$ and
have side $1/m_\la\sim r_\la/(2d)$, assuming $r_\la\to0$,
which holds by (\ref{fromThetacond}) if $p_\la= \Omega(\la^{-\eps
})$ for
some $\eps\in(0,1)$.

Given $\la$,
for $z \in\Lambda_\la$ let us
say the cube $Q_z$ is \emph{blue} if:
(i) $\Po_\la\cap Q_z \neq\varnothing$ and
(ii) all vertices of $\Po_\la\cap B(m_\la^{-1}z;r_\la/\eta)$
lie in the same connected component of
$G_{\phi_\la}(\Po_\la\cap B(m_\la^{-1} z; 2 r_\la/\eta))$.
If a cube is not blue, let us say it is \emph{green}.
If $Q_z$ is blue (resp., green), we shall also say
$\overline{Q}_z$ and $z$
are blue (resp., green). More prosaically we shall put
$Y_{\la,z} =1$ if $z$ is blue
and
$Y_{\la,z} =0$ if $z$ is green.

\begin{lemm}\label{lemblue}
Suppose
$p_\la= \Omega(\la^{-\eps})$ with $0 < \eps< (9d)^{-d} \eta/
K(\eta)$.
Then
\[
\sup_{z \in\Lambda_{\lambda}} \Pr[ Y_{\la,z} =0 ] = O\bigl(
\la^{-\eps}\bigr).
\]
\end{lemm}

\begin{pf}
First note that
$\card(\Po_\la\cap Q_z)$
is Poisson with mean $\la/ m_\la^d \sim(2d)^{-d} \la r_\la^d
\geq(2d)^{-d} \la I(\phi_\la)/K(\eta)
$, where the inequality comes from
(\ref{0729b}). Hence by (\ref{prcond0}) the probability that
condition (i) (in the definition of blue)
fails is $O(\la^{-(3d)^{-d}/K(\eta)})$, uniformly
over $z \in\Lambda_\la$.
We need a similar bound for
the probability that
condition (ii) fails.

Let $\xi_\la$ be Poisson with parameter
$2 \la/ m_\la^d$. We claim that the
{Erd\H{o}s--R\'enyi }graph $G(\xi_\la,\eta p_\la)$ satisfies
%
\begin{equation}
\Pr\bigl[ G(\xi_\la,\eta p_\la) \notin\K\bigr] = O\bigl(
\la^{-\eps}\bigr). \label{0717a}
\end{equation}
Indeed, by the Mecke formula followed by (\ref{fromThetacond}),
(\ref{0729b})
and (\ref{prcond0}), the expected number of isolated vertices in
$G(\xi_\la,\eta p_\la)$
is given by
\begin{eqnarray*}
O \bigl( \la r_\la^d \exp\bigl( - (3d)^{-d}
\la r_\la^d \eta p_\la\bigr) \bigr) &=& O\bigl(
\la^{2 \eps} \exp\bigl( - (3d)^{-d} \eta\la I(
\phi_\la) /K(\eta) \bigr) \bigr)
\\
&=& O 
\bigl( \la^{2 \eps- (3d)^{-d}\eta/K(\eta)} \bigr), 
\end{eqnarray*}
which is $O(\la^{-\eps})$
by the condition on $\eps$.
Thus the probability that
$G(\xi_\la,p_\la)$
has an isolated vertex is $O(\la^{-\eps})$,
and by the proof of \cite{Boll}, Theorem~7.2,
we have (\ref{0717a}).
Hence, for each
pair of neighboring sites
$z',z'' \in\Lambda_\la$,
the graph $G_{\phi_\la}(\Po_\la\cap(Q_{z'} \cup Q_{z''}))$ is
connected with probability $1-O(\la^{-\eps})$.
Condition (ii) holds if
$G_{\phi_\la}(\Po_\la\cap(Q_{z'} \cup Q_{z''}))$ is
connected for each pair of neighboring sites
$z',z''$ lying in $B(z, 2m_\la r_\la/\eta) \cap\Lambda_\la$,
and the number of such pairs is bounded independently of $z$ and
$\la$.
Therefore by
the union bound, condition
(ii) holds
with probability $1-O(\la^{-\eps})$,
as claimed.
\end{pf}

We say a set $S \subset\Lambda_\la$ is $*$-connected if
for any $x,y \in S$,
there is a path $(x_0,x_1,\ldots,x_k)$
with $x_0=x$, $x_k=y$ and $x_i \in
S
$ and
$\|x_i -x_{i-1}\|_\infty= 1$
for $1 \leq i \leq k$ (so diagonal steps in the path are allowed).
For bounded nonempty $U \subset\R^d$, we define the
\emph{$\ell_\infty$-diameter} of $U$ to be $\sup_{x,y\in U} \|y-x\|
_\infty$.
Given $\la,\rho>0$, let $H_\la^\rho$ be the event that there
is a $*$-connected set of green sites in $\Lambda_\la$
of $\ell_\infty$-diameter at least
$\rho$.

\begin{lemm}
\label{lemsea}
Suppose
for some $\eps\in(0,(9d)^d \eta/K(\eta))$ that
$p_\la=\Omega(\la^{-\eps})$.
Then
there exists $\rho>0$ such that
$\Pr[H_\la^\rho]\to0$ as $\la\to\infty$.
\end{lemm}

\begin{pf}
For $\la>0, n \in\N$, let
${\cal T}_{\la,n}$ denote the set of $*$-connected sets
$\gamma\subset\Lambda_\la$ with
$n$ elements.
Then there exists a constant $A$ such that for all $\la$ and $n$,
we have
$\card( {\cal T}_{\la,n} ) \leq m_\la^d A^n$; see, for example,
\cite{Penbk}, Lemma~9.3.
Also $r_\la^{-d} = \Theta(\la p_\la/\log\la)$ by (\ref{fromThetacond}),
and hence there exists $\la_0 \in(0,\infty)$ such that
for $\la\geq\la_0$ we have
that $m_\la^{d} \leq\la$
so that
$\card( {\cal T}_{\la,n} ) \leq\la A^n$
for all $n \in\N$.

The random field $(Y_{\la,z},z \in\Lambda_\la)$
has finite range dependency: there exists $\la_1 \in[ \la_0,\infty)$
such that
the range may be taken to be $11d/\eta$, for all $\la\geq\la_1$.
For example, if $|z-z'| \geq11 d/\eta$,
then $|m_{\la}^{-1}z - m_{\la}^{-1}z'| \geq5 r_\la/ \eta$,
and therefore $Y_{\la,z}$ is independent of
$Y_{\la,z'}$.
Therefore there is a constant $M:= M(d,\eta)$
such that for any $\la\geq\la_1$ and any $S \subset\Lambda_\la$,
we can find $S' \subset S$ with
$\card(S') = \lceil\card(S)/M\rceil$, such that
the variables $(Y_{\la,z})_{z \in S'}$ are
mutually independent.
Hence by Lemma~\ref{lemblue} there is a further
constant $C$ such that for all such $S$ we have
\[
P\biggl[ \bigcap_{z \in S} \{ Y_{\la,z} =0\}
\biggr] \leq\bigl(C \la^{-\eps}\bigr)^{(\card S) /M}.
\]

Let $\rho\in\N$.
If $H^\rho_\la$ occurs, then
there exists $S \in{\cal T}_{\lambda,\rho}$ such that
$ Y_{\la,z} =0 $ for all $z \in S$.
Hence for $\rho\in\N$ and
$\la\geq\max(\la_1,(CA^M)^{2 /\eps})$, we have
\[
\Pr\bigl[H_\la^\rho\bigr] \leq \Pr\biggl[ \bigcup
_{S \in{\cal T}_{\la,\rho}} \bigcap_{z \in S} \{
Y_{\la,z} =0\} \biggr] \leq\la A^{ \rho} \bigl(C
\la^{-\eps}\bigr)^{\rho/M} \leq\la^{1 -\eps\rho/(2M)}.
\]
Taking
$\rho> 2 M/\eps$, we have the result.
\end{pf}

Given disjoint nonempty connected subsets $U$ and $V$ of
$\Gamma$, we define the \emph{exterior boundary}
of $U$ relative to $V$ as follows. Let
$V'$ be the connected component of
$\Gamma\setminus U$ that contains $V$,
and let $U':= \Gamma\setminus V'$. Loosely speaking,
$U'$ is obtained from $U$ by filling in all the holes
in $U$, except the one containing $V$.
Define the exterior boundary of $U$ relative
to $V$ to be the intersection of the
closure of $U'$ with that of $V'$.

The exterior boundary of $U$ relative to
$V$ is a subset of the boundary of $U$. Moreover
it is a connected set,
by a unicoherence argument (see \cite{Penbk}),
because the closures of $U'$ and $V'$ are connected
sets whose union is $\Gamma$.

We claim that for $ 0 <a <1$, if both $U$ and $V$
have $\ell_\infty$-diameter greater than $a$,
then so does the exterior boundary of $U$ relative to
$V$. Indeed, if not, then there exists a rectilinear
cube ${\mathbf C}$
of side $a$ that contains the exterior
boundary of $U$ relative to $V $, but then we could
pick $u \in U \setminus{\mathbf C}$ and $v \in V \setminus{\mathbf C}$,
and a continuous path
from $u$ to $v$ in $\Gamma$ avoiding ${\mathbf C}$. Somewhere
on this path would lie a point in the exterior boundary
of $U$ relative to $V$, a contradiction.

\begin{lemm}
\label{lembigblue}
Let $\la>0$, $\rho\in\N$ with $ \rho< m_\la$,
and suppose $H_\la^\rho$ does not occur. Then
there exists a
$*$-connected
component of the set of blue sites in $\Lambda_\la$ of
$\ell_\infty$-diameter $m_\la-1$.
This component is
unique, and there
is no other $*$-connected
component of the set of blue sites in $\Lambda_\la$ of
$\ell_\infty$-diameter
$\rho$ or more.
\end{lemm}

\begin{pf}
Let $B_\la$ denote the union of all
the cubes $Q_z,z \in\Lambda_\la$ that are blue, and
let $G_\la$ denote the union of all
the cubes $Q_z,z \in\Lambda_\la$ that are green.
Let $U$ be the component of $G_\la\cup(\{0\} \times[0,1]^{d-1})$
that contains $\{0\} \times[0,1]^{d-1}$, and let $V$ be the
component of $B_\la\cup( \{1\} \times[0,1]^{d-1})$
that contains $\{1\} \times[0,1]^{d-1}$.
Then $U$ and $V$ are disjoint connected subsets of $\Gamma$.
Assuming $H_\la^\rho$ does not occur,
$U$ does not extend to
$\{1\} \times[0,1]^{d-1}$. Hence
the union of blue cubes $\overline{Q}_z$
having nonempty intersection with the
exterior boundary of $U$ relative to $V$ is
connected and has $\ell_\infty$-diameter
$1$, and the first assertion (existence) in the statement
of the lemma follows.

Suppose there were two
$*$-connected components
of the set of blue sites
of $\ell_\infty$-diameter at least $\rho$, denoted
$U$ and $V$, say.
Let $U^*$ be the union of
the cubes $\overline{Q}_z,z \in U$, and define $V^*$ similarly.
Then $U^*$ and $V^*$ are connected disjoint regions of $\Gamma$, of
$\ell_\infty$-diameter
at least
$ (\rho+1) /m_\la$.
The union of
green cubes $\overline{Q}_y$ having nonempty intersection
with the exterior boundary of $U^*$ relative to
$V^*$ would be a connected region of $\ell_\infty$-diameter at
least
$ (\rho+1)/ m_\la$,
and the corresponding
set of sites in $\Lambda_\lambda$ would be a $*$-connected set of
green sites of diameter at least
$ \rho$,
contradicting the assumed nonoccurrence of event $H^\rho_\la$.
This demonstrates the second assertion (uniqueness) in
the statement of the lemma.
\end{pf}

We shall refer to the unique $*$-connected
blue component of $\ell_\infty$-diameter \mbox{$m_\la-1$},
identified in Lemma~\ref{lembigblue}, as the \emph{sea}.
All vertices of $\Po_\la$ lying in cubes $Q_z$ with $z$ in the sea lie
in the same component of $G_{\phi_\la}(\Po_\la)$,
which we call the \emph{sea-component}.

Given $\la>0$, $\rho>0$, define the event
\begin{eqnarray*}
&& F_\la^\rho= \bigl\{ \exists x,y \in\Po_\la
\dvtx \min\bigl( D_{\phi_\la
}(x,\Po_\la), D_{\phi_\la}(y,
\Po_\la)\bigr) > \rho,
\\
&&\hspace*{127pt} C_{\phi_\la}(x,\Po_\la) \neq C_{\phi_\la}(y,
\Po_\la) \bigr\}.
\end{eqnarray*}

\begin{lemm}
\label{lemmid}
Let $0 < \eps< (9d)^{-d} \eta/K(\eta)$.
There exists a constant $\rho\in\N$, such that if for
some $\alpha>0$, we have (\ref{prcond0}) and
also $p_\la= \Omega(\la^{-\eps})$, then
$
P[
F_\la^{\rho r_\la}
] \to0
$ as $\la\to\infty$.
\end{lemm}

\begin{pf}
Let $\rho\in\N$.
Suppose that
$F_\la^{\rho r_\la}$ occurs and
$H_\la^\rho$ does not. Then there exists $U \subset\Po_\la$
such that $U$ is the vertex-set of a component
of
$G_{\phi_\la}(\Po_\la)$ that is disjoint
from the sea-component, but has diameter greater than $\rho r_\la$,
and hence has $\ell_\infty$-diameter greater than
$\rho r_\la/\sqrt{d}$.

Let $\tilde{U}$
denote the union of closed Euclidean balls
of radius
$r_\la/(2 \eta)$ centered on the
vertices of
$U$.
This is a connected subset of $\R^d$
because $\rho_0 (\phi_\la) \leq\eta^{-1}r_\la$ by
(\ref{eqPhi0}), and
therefore for each pair of vertices $y,y'$ connected
by an edge of $G_{\phi_\la}(\Po_\la)$,
we have $|y-y'| \leq r_\la/\eta$.
Also $\tilde{U}$
has $\ell_\infty$-diameter of at least $\rho r_\la/\sqrt{d}$.

We claim there is no $x \in U$ and $z $ in the sea
such that $|x - m_\la^{-1}z | \leq\eta^{-1} r_\la$.
For if there were such a pair, then by the
definition of blue, $x$ would lie in the same component
as the vertices of $\Po_\la$ in $Q_z$,
so $U$ would be part of the sea-component,
a contradiction.

Let $S$ be the union of cubes $\overline{Q}_z$
with $z$ in the sea.
The set $S$ is connected, and disjoint from
$\tilde{U}$ by the preceding claim, since
the cubes have diameter at most $r_\la/(2 \sqrt{d})$;
let $\partial_\mathrm{ ext} \tilde{U}$
denote the exterior boundary of $
\tilde{U}$ relative to $S$. This has
$\ell_\infty$-diameter at least $\rho r_\la/\sqrt{d}$.

Now let $\Delta_\mathrm{ ext}\tilde{U}$ be the set
of sites $z \in\Lambda_\la$ such that the corresponding
cubes $\overline{Q}_z$ have nonempty intersection with
$\partial_\mathrm{ ext} \tilde{U}$.
Since $\partial_\mathrm{ ext} \tilde{U}$ is connected,
the set $\Delta_\mathrm{ ext}\tilde{U}$ is $*$-connected.
Also
$ \card( \Delta_\mathrm{ ext}\tilde{U}) \geq(\rho r_\la/\sqrt
{d})m_\la-1
\geq\rho$.

We claim that none of the squares $Q_z, z \in
\Delta_\mathrm{ ext}\tilde{U}
$,
is blue.
This is because by definition,
each such $Q_z$ intersects with
$\partial_\mathrm{ ext}\tilde{U}$,
and therefore lies
at a distance of at most $r_\la/(2 \eta)$ from some vertex of
$U$
(at $X$, say).
Then by the triangle inequality $|X-m_\la^{-1} z | \leq
r_\la/(2\eta) + r_\la/(2 \sqrt d) \leq r_\la/\eta$, and
if $Q_z$ were blue,
it would contain at least one vertex of $\Po_\la$,
and this would be in the same component of
$G_{\phi_\la}(\Po_\la)$ as all the vertices within
distance $r_\la/\eta$ of $m_\la^{-1}z$, including $X$. Hence
$Q_z$ would include a vertex of
$U$,
but then it would be contained in the interior of
$\tilde{U}$, and so would
have empty intersection with
$\partial_\mathrm{ ext} \tilde{U}$,
a contradiction.

Thus $\Delta_\mathrm{ ext}\tilde{U}$ is
a $*$-connected set of cardinality at least $\rho$,
all of whose elements are green. This contradicts
the assumed nonoccurrence of $H_\la^\rho$.
Thus $F_\la^{\rho r_\la} \subset H_\la^\rho$,
and the result follows from Lemma~\ref{lemsea}.
\end{pf}

\begin{pf*}{Proof of Theorem~\ref{prop0711c}}
Set $\eps= \frac{1}{2} \min((9d)^{-d}\eta/K(\eta),2^{-d -3}) $.
Given $\rho>0$, if $ L_2 (G_{\phi_\la}(\Po_\la)) > 1 $,
then either $E_\la^{\rho r_\la}$ or $F_\la^{\rho r_\la}$
occurs.
If $p_\la= \Omega(\la^{-\eps})$,
result~(\ref{0711f}) follows from
Proposition~\ref{lemd3c} and Lemma~\ref{lemmid}.
If $p_\la= O(\la^{-\eps})$,
(\ref{0711f}) follows from Proposition~\ref{Poconnthm}.
\end{pf*}

\section{De-Poissonization}
\label{secdepo}

In this section we shall
complete the proof of Theorems
\ref{thm2a}, \ref{thm1} and
\ref{thm2}. We start with
the case $\alpha\in(0,\infty)$ of Theorem~\ref{thm2}.
All integrals in this section are over $\Gamma$ unless
specified otherwise.

\begin{prop}
\label{theodepo}
Suppose $\alpha\in(0,\infty)$ and $(\phi_n)$ satisfy (\ref{prcond1})
as $n \to\infty$ along some subsequence of $\N$,
and for some $\eta\in(0,1]$ we have $\phi_n \in\Phi_{d,\eta}$ for
all $n$.
Then for $k \in\N_0$, (\ref{eq1b}) holds as $n \to\infty$
along the same subsequence.

If also
$\phi_n \in\Phi_{d,\eta}^0$ for all $n$,
then along the same subsequence we have
%
\begin{equation}
\label{main2} \lim_{n \to\infty} \Pr\bigl[L_2
\bigl(G_{\phi_n}(\X_n) \bigr) \leq1 \bigr] =1.
\end{equation}
\end{prop}
\begin{pf}
Let $\la(n) = n- n^{3/4}$ and $\mu(n):= n+ n^{3/4}$.
Let $\Po_{\la(n)}, \X_n, \Po_{\mu(n)}$ be coupled as
follows. Let $X_1,X_2,\ldots $ be a sequence
of independent random vectors uniformly
distributed over $\Gamma$.
Independently,
let $Z$ and $Z'$ be Poisson distributed
random variables
with parameter
$\la(n) $ and $\mu(n)-\la(n)$, respectively, independently of
each other and of $(X_1,X_2,\ldots)$;
set $ \Po_{\la(n)}:= \{X_1,\ldots,X_{Z}\}$,
and set $ \Po_{\mu(n)}:= \{X_1,\ldots,X_{Z+Z'}\}$
and $\X_n:= \{X_1,\ldots,X_{n}\}$.
By Chebyshev's inequality,
w.h.p. $\Po_{\la(n)} \subset\X_n \subset\Po_{\mu(n)}$.

Without loss of generality, assume $\rho_\eta(\phi_n) \leq\sqrt{d}$.
By (\ref{Thetacond}),
\[
\exp \biggl( n^{3/4} \int_\Gamma
\phi_n(y-x) \,dy \biggr) = \exp\bigl( n^{-1/4} \times\Theta(
\log n) \bigr) =1+o(1),
\]
uniformly over $x \in\Gamma$,
and therefore
the sequence $(\phi_n)_{n \in\N}$ satisfies
%
\begin{equation}
\la(n) \int_\Gamma\exp \biggl( - \la(n) \int
_\Gamma\phi_n (y-x) \,dy \biggr) \,dx \to\alpha.
\label{prminus}
\end{equation}

Let $A_n$ be the union of the event that
at least one of the added vertices of $\Po_{\mu(n)} \setminus
\Po_{\la(n)}$
is not connected to any of the vertices of
$ \Po_{\la(n)}$, and the event that
at least one of the added vertices of $\Po_{\mu(n)} \setminus
\Po_{\la(n)}$
is connected to one of the
isolated vertices of
$G_{\phi_n}(\Po_{\la(n)}) $.

By the Mecke equation,
the expected number of added vertices
that are isolated from
all the vertices of $\Po_{\la(n)}$
equals
$ 2 n^{3/4} \int
\exp (- \lambda(n) \int\phi_n(y-x)
\,dy  ) \,dx $,
which tends to zero
by (\ref{prminus}).
Also, the expected number of isolated vertices in $G_{\phi_n}(\Po
_{\la(n)})$
that are connected to at least one of the
added vertices is bounded by
\[
\bigl(n - n^{3/4}\bigr) \int_\Gamma\exp \biggl(-
\bigl(n - n^{3/4}\bigr) \int_\Gamma
\phi_n (y-x) \,dy \biggr) \times2 n^{3/4} I(
\phi_n) \,dx,
\]
and by (\ref{prminus}) and (\ref{Thetacond}) this
tends to zero.
Hence $\Pr[A_n] = o(1)$.
By Theorem~\ref{Poissecthm}
we have that
\[
\Pr\bigl[N_0\bigl( G_{\phi_{n}}(\Po_{\la(n)})\bigr) =k
\bigr] \to e^{-\alpha} \alpha^k/k!,\qquad k \in\N_0.
\]
Also
$\Pr[ N_0(G_{\phi_n}(\X_n)) \neq
N_0(G_{\phi_n}(\Po_{\la(n)})) ] \leq\Pr[ A_n] +
\Pr[ \{Z \leq n \leq Z + Z'\}^c]$, which
tends to 0, and
(\ref{eq1b}) follows.

Now suppose
$\phi_n \in\Phi_{d,\eta}^0$ for all $n$.
If $L_2(G_{\phi_n}(\X_n)) >1 $,
then either $Z > n$, or $Z + Z' < n$,
$L_2(G_{\phi_n}(\Po_{\la(n)})) > 1$, or $A_n$ occurs.
By
Theorem~\ref{prop0711c}, all of these events have
vanishing probability, and (\ref{main2}) follows.
\end{pf}

Next we consider the case with $\alpha\in\{0,\infty\}$.

\begin{prop}
\label{prop6inf}
Suppose $\alpha\in\{0,\infty\}$, $\eta\in(0,1]$
and $(\phi_n)$ satisfy (\ref{prcond1}) as $n \to\infty$ along
some subsequence of $\N$,
and $\phi_n \in\Phi_{d,\eta}$
for all $n$.
If $\alpha=0$, then $\Pr[N_0(G_{\phi_n} (\X_n)) =0] \to1$,
and
if $\alpha=\infty$, then for all
$k \in\N_0$,\break  $\Pr[N_0(G_{\phi_n} (\X_n)) =k] \to0$,
as $n \to\infty$ along the same subsequence.
\label{prdepo2}
\end{prop}

\begin{pf} (i)
Let $I_n(\phi_n)$ denote the left-hand side of (\ref{prcond1}).
Then
\begin{eqnarray*}
\E N_0 \bigl(G_{\phi_n} (\X_n)\bigr) &=& n \int \,dx
\biggl(1 - \int\phi_n(y-x) \,dy \biggr)^{n-1}
\\
&\leq& n \int \,dx \biggl( \exp \biggl( - (n-1) \int\phi_n(y-x) \,dy
\biggr) \biggr) \leq e I_n(\phi_n).
\end{eqnarray*}
Therefore, by Markov's inequality, if
$\alpha=0$, we have $
\Pr[ N_0 (G_{\phi_n} (\X_n) ) \geq1] \to0$.

Now suppose $\alpha=\infty$.
We seek to interpolate a ``larger'' connection function
than $\phi_n$ that is still
in $\Phi_{d,\eta}$.
For $s >1$ and $\phi\in\Phi_{d,\eta}$,
define $\phi^{(s)}$ as follows.
Let $s_0(\phi) = 1/ \mu(\phi)$.
For $1 \leq s \leq s_0(\phi)$,
set $\phi^{(s)} (x):= s \phi(x)$, for $x \in\R^d$.
Note $\mu(\phi^{(s_0(\phi))})=1$.
For $s \geq s_0(\phi)$, define
\[
\phi^{(s)}(x): = \cases{ 1, &\qquad $\mbox{if } |x| < s-s_0(\phi),$
\vspace*{2pt}
\cr
\phi_{s_0}(x), &\qquad  $\mbox{if } |x| \geq
s-s_0(\phi) $. }
\]
Let $s_1(\phi):= \sqrt{d} + s_0(\phi)$.
If $\phi\in\Phi_{d,\eta}$, then for each $s \in[1, s_1(\phi)]$ the
connection function $\phi^{(s)}$ is also in $\Phi_{d,\eta}$.

For each $n \in\N$ define the function
\[
\tilde{f}_n(s): = n \int\exp \biggl( - n \int\phi_n^{(s)}
(y-x) \,dy \biggr) \,dx,
\]
which is continuous and nonincreasing on $1 \leq s \leq s_1(\phi_n)$.
By assumption $\tilde{f}_n(1) \to\infty$ as $n \to\infty$, while
$\tilde{f}_n(s_1(\phi_n)) = n e^{-n}$. Therefore
by the intermediate
value theorem,
given any finite $\beta>0$, for large enough $n$ we can pick
$s(n) \in[1,s_1(\phi_n)]$ with $\tilde{f}_n(s(n)) =\beta$.
Then by Proposition~\ref{theodepo}, for $k \in\N_0$ we have
\[
\Pr\bigl[N_0\bigl(G_{ \phi_n^{(s(n))}}(\X_n)\bigr) \leq k
\bigr] \to e^{-\beta} \sum_{j=0}^k
\beta^j/j!.
\]
By an obvious coupling,
$\Pr[N_0(G_{\phi_n^{(s)}}(\X_n)) \leq k] $
is nondecreasing in $s$, and therefore
since $\beta>0$ is arbitrary, we have
$\Pr[N_0(G_{ \phi_n}(\X_n)) \leq k] \to0$.
\end{pf}

\begin{pf*}{Proof of Theorem~\ref{thm1}}
Let $\eta\in( 0,1]$. To prove (\ref{eqmain}),
it suffices
to prove that for any sequence $(\phi_n)_{n \in\N}$ of
connection functions in $\Phi_{d,\eta}^0$,
we have
%
\begin{equation}
\lim_{n \to\infty} \Pr\bigl[ \bigl\{N_0
\bigl(G_{\phi_n} (\X_n)\bigr) =0\bigr\} \setminus \bigl
\{G_{\phi_n}(\X_n) \in\K\bigr\} \bigr] =0. \label{wkmain}
\end{equation}
Define
$I_n:= I_n(\phi_n):= n \int\exp ( - \int\phi_n(y-x) \,dy
 ) \,dx $.
Consider the three cases where:
(i) $I_n$ tends to a finite limit as $n \to\infty$ along
some infinite subsequence of $\N$;
(ii) $I_n \to\infty$ as $n \to\infty$ along some infinite
subsequence of $\N$;
(iii)
$I_n \to0$ as $n \to\infty$ along some infinite subsequence of $\N$.
At least one of these cases holds, and it suffices
to show that in each case (\ref{wkmain}) holds along the same subsequence.

In case (i), we have (\ref{wkmain}) at once because
of (\ref{main2}).
In case (ii), with $I_n \to\infty$,
by Proposition~\ref{prdepo2}
we have $\Pr[N_0(G_{\phi_n}(\X_n)) =0] \to0$, and hence
(\ref{wkmain}) holds.

Consider case (iii) with $I_n \to0$ along a subsequence.
For $n \in\N$,
define
\[
f_n(a):= n \int\exp \biggl(-a n \int\phi_n(y-x) \,dy
\biggr) \,dx,
\]
which is a continuous and nonincreasing function on $0 \leq a \leq1$.
For each $a \in[0,1]$ the connection function $a \phi_n$ is in $\Phi
_{d,\eta}$.

By assumption $f_n(1) \to0$ as $n \to\infty$,
while $f_n(0)=n$.
Therefore given $\eps>0$, by the intermediate value theorem,
for all large enough $n$ in the subsequence
we can choose $a_n \in[0,1]$ such
that $f_n(a_n) = \eps$. Then by Proposition~\ref{theodepo}
we have
\[
\Pr\bigl[N_0\bigl(G_{a_n \phi_n}(\X_n)\bigr) =0\bigr]
\to e^{-\eps};\qquad \Pr\bigl[G_{a_n \phi_n}(\X_n) \in\K\bigr]
\to e^{-\eps}.
\]
By an obvious coupling,
$\Pr[G_{a \phi_n}(\X_n) \in\K] $
is nondecreasing in $a$, and therefore
since $\eps$ is arbitrary, we have
$\Pr[G_{\phi_n}(\X_n) \in\K] \to1$,
so (\ref{wkmain}) holds.
\end{pf*}

\begin{pf*}{Proof of Theorem~\ref{thm2}} Equation (\ref{eq1b}) follows from
Proposition~\ref{theodepo}, and the next sentence follows
from Proposition~\ref{prop6inf}.
Then (\ref{eq1c}) follows from Theorem~\ref{thm1}.
\end{pf*}

\begin{pf*}{Proof of Theorem~\ref{thm2a}}
The result follows from Theorem~\ref{thm2}.
\end{pf*}

\section{Equivalence of thresholds}
\label{secequiv}

In this section we prove Theorem~\ref{threshthm}; that is, we prove
that for any $[0,1]$-valued sequence $(p_n)_{n \in\N}$
with $p_n = \omega( (\log n)/n)$,
we have
\[
\lim_{n \to\infty} P\bigl[\tau_n(p_n) =
\sigma_n(p_n)\bigr] =1,
\]
where for $p \in[0,1]$,
as described in Section~\ref{secmainres}
we set
\[
\tau_n(p):= \inf\bigl\{r\dvtx G_{r,p}(\X_n)
\in\K\bigr\};\qquad \sigma_n(p):= \inf\bigl\{r\dvtx N_0\bigl(
G_{r,p}(\X_n)\bigr) = 0 \bigr\}.
\]

Clearly $\sigma_n(p_n) \leq\tau_n(p_n)$,
so we need to show that $ \Pr[\sigma_n(p_n) < \tau_n(p_n)]$
tends to zero.
Given $p_n$ and given $\alpha>0$, define $r_n(\alpha)$ by
$I_n(\phi_{r_n(\alpha),p_n} ) = e^{- \alpha}$, where
$I_n(\phi):= n \int_\Gamma\exp( - n \int_\Gamma\phi(y-x) \,dy)\,dx$.
For each $\alpha$ we have from (\ref{eqmain})
that
%
\begin{equation}
P\bigl[ \sigma_n(p_n) \leq r_n(\alpha) <
\tau_n(p_n) \bigr] \to0. \label{0109a}
\end{equation}
Note that $r_n(\alpha)$ is nondecreasing in $\alpha$.
Let $\alpha< \beta$.
Suppose
\[
r_n (\alpha) < \sigma_n(p_n) <
\tau_n(p_n) \leq r_n(\beta).
\]
Assume the inter-point distances are all distinct.
Consider adding the edges of $G_{\sqrt{d},p} (\X_n)$
one by one (starting from the graph with no edges)
in order of increasing Euclidean length.

Then precisely one pair of points of $\X_n$,
say $X$ and $Y$, satisfies $|X-Y| = \tau_n(p_n) $,
and by the definition of $\tau_n(p_n)$, $X$ and $Y$ lie
in different components
just before adding the edge
between them, so they lie in different components
of $G_{r_n(\alpha),p_n}(\X_n)$.
Assuming $L_2(G_{r_n(\alpha),p_n} (\X_n)) \leq1$
[which has high probability by (\ref{main2})],
either $X$ or $Y$ (say $X$) is isolated in
$G_{r_n(\alpha),p_n}(\X_n)$, but
$X$ is nonisolated
in
$G_{\sigma_n(p_n),p_n}(\X_n)$ by the definition of $\sigma_n(p_n)$.
Therefore since we are assuming $\tau_n(p_n) \leq r_n(\beta)$,
we have that $X$ is connected to at least two
points of $\X_n$, at distances between
$r_n(\alpha)$ and $r_n(\beta)$.
Thus $N_{\alpha,\beta}(n) >0$, where
$N_{\alpha,\beta}(n)$ denotes
the number of vertices of $\X_n$
having no incident edge in
$G_{\sqrt{d},p_n}(\X_n)$ of (Euclidean) length at most
$ r_n(\alpha)$
but at least
two incident edges of length at most $r_n(\beta)$.

Let $\la(n)$ and $\mu(n)$, and
the coupling of $\Po_{\la(n)},\X_n$, and $\Po_{\mu(n)}$
be as in the preceding section.
Let $N'_{\alpha,\beta}(n)$
be the number of vertices of $\Po_{\mu(n)}$
having no incident edge [in $G_{\sqrt{d},p_n}(\Po_{\mu(n)})$]
of length at most $ r_n(\alpha)$
with the other endpoint in $\Po_{\la(n)}$,
but at least
two incident edges
of length at most $ r_n(\beta)$
(with the other endpoint in $\Po_{\mu(n)}$).
If $\Po_{\la(n)} \subset\X_n \subset\Po_{\mu(n)}$ (which happens
with high probability), then $N'_{\alpha,\beta}(n) \geq N_{\alpha
,\beta}(n)$.
Thus
%
\begin{equation}\qquad
\limsup_{n \to\infty} \Pr\bigl[ r_n(\alpha) <
\sigma_n(p_n) < \tau_n(p_n) \leq
r_n(\beta) \bigr] \leq\limsup_{n \to\infty} \Pr
\bigl[N'_{\alpha,\beta
}(n) >0\bigr]. \label{1104a}
\end{equation}

With $|\cdot|$ denoting Lebesgue measure, by the Mecke formula
we have
\[
\E\bigl[N'_{\alpha,\beta}\bigr] = \bigl(n+n^{3/4}\bigr)
\int_{\Gamma} e^{-\lambda(n) p_n |B(x;r_n(\alpha)) \cap\Gamma|} \times \bigl(1 -
e^{-w_n(x)}\bigl(1 + w_n(x)\bigr)\bigr) \,dx,
\]
where $w_n(x)$ denotes the mean number of edges
of length in the range $(r_n(\alpha), r_n(\beta)]$ incident
to a point at $x$.
Now,
$e^{w} -1 -w \leq w^2 e^{w}$ for any $w \geq0$.
Hence
%
\begin{equation}
\E\bigl[N'_{\alpha,\beta}\bigr] \leq \bigl(n+n^{3/4}
\bigr) \int e^{-\lambda(n)p_n|B(x;r_n(\alpha)) \cap\Gamma|} \times w_n(x)^2 \,dx.
\label{sq1}
\end{equation}
By (\ref{fromThetacond}) and the condition $p_n = \omega((\log n)/n)$,
we have $r_n(\beta) \to0$.
Writing $V_\alpha(x)$ for $|B(x;r_n(\alpha)) \cap\Gamma|$
we have
\begin{eqnarray*}
e^{-\alpha}& = &\lim_{n \to\infty} \biggl( n \int_\Gamma
\exp\bigl(-n p_n V_\beta(x) + np_n
\bigl(V_\beta(x) - V_\alpha(x) \bigr)\bigr) \,dx \biggr)
\\
&\geq&\limsup_{n \to\infty} \biggl( n \int_\Gamma
\exp \bigl( - n p_n V_\beta(x) + n p_n
\pi_d \bigl(r_n(\beta)^d -r_n(
\alpha)^d\bigr)/2^d \bigr) \,dx \biggr)
\\
&=& e^{-\beta} \exp\Bigl( \limsup_{n \to\infty} \bigl[n
p_n \pi_d \bigl(r_n(\beta)^d
-r_n(\alpha)^d\bigr)/2^d\bigr]\Bigr)
\end{eqnarray*}
so that
\[
\limsup_{n \to\infty} n p_n \bigl(r_n(
\beta)^d -r_n(\alpha)^d \bigr)
\leq2^d (\beta-\alpha)/\pi_d.
\]
Therefore, since
\[
w_n(x) \leq\mu(n) p_n \pi_d
\bigl(r_n(\beta)^d - r_n(\alpha)^d
\bigr),
\]
we have $\limsup_{n \to\infty} \sup_{x \in\Gamma}
w_n(x) \leq2^d(\beta- \alpha)$,
so that by (\ref{sq1}) and a similar argument to
(\ref{prminus}),
$
\limsup_{n \to\infty}
\E[N'_{\alpha,\beta}] \leq
2^{2d} (\beta-\alpha)^2 e^{-\alpha}$,
so that by (\ref{1104a}),
%
\begin{equation}
\limsup_{\la\to\infty} P\bigl[r_n (\alpha) <
\sigma_n(p_n) < \tau_n(p_n) \leq
r_n(\beta)\bigr] \leq2^{2d} (\beta- \alpha)^2
e^{- \alpha}. \label{0109b}
\end{equation}
Now we argue as in
\cite{kcon}, pages 163--164
or \cite{Penbk}, pages 304--305. Let $\eps>0$. Choose $\alpha_0 <
\alpha_1
< \cdots< \alpha_I$ such that
$\exp(- e^{-\alpha_0} ) < \eps$, and $1 -
\exp(- e^{-\alpha_I} ) < \eps$, and
also
\[
2^{2d} \sum_{i=1}^I
\bigl(r_n(\alpha_i) - r_n(
\alpha_{i-1}) \bigr)^2 e^{-\alpha_{i-1} } < \eps.
\]
Then by the union bound,
\begin{eqnarray*}
P [ \sigma_n < \tau_n ] &\leq& P\bigl[\sigma_n
\leq r_n (\alpha_0) \bigr] + P\bigl[\sigma_n
> r_n (\alpha_I) \bigr]
\\
&&{}+ \sum_{i=1}^I \bigl( P \bigl[
\sigma_n \leq r_n(\alpha_i) <
\tau_n\bigr] + P \bigl[ r_n(\alpha_{i-1}) <
\sigma_n < \tau_n \leq r_n(
\alpha_{i})\bigr] \bigr).
\end{eqnarray*}
Since $\sigma_n \leq r$ if and only if $N_0(G(\X_n,r) ) =0$,
it follows from (\ref{eq1b}) of Theorem~\ref{thm2},
along with (\ref{0109a}) and (\ref{0109b}), that
$
\limsup_{n \to\infty}
P [ \sigma_n < \tau_n ] \leq3 \eps$,
and since $\eps>0$ is arbitrary, this completes the proof.

\section{The choice of \texorpdfstring{$\phi$}{phi}}
\label{secchoice}

In this section, we prove Theorem~\ref{thconds} (among other things).
That is, we identify conditions for a sequence of connection
functions $\phi_n$ to satisfy (\ref{prcond1}) for some $\alpha\in
(0,\infty)$.
We consider only the case with $d=2$ and
$\phi_n \in\Phi_{2,\eta} \cap\Psi_2$
for some $\eta\in(0,1]$, where
$\Psi_2$ is defined by (\ref{0905a}).

Assume $d=2$. Fix $\eta>0$, and
choose $\phi_n \in\Phi_{2,\eta} \cap\Psi_2 $
for each $n >0$.
Set
\[
r_n:= \rho_\eta(\phi_n);\qquad p_n:=
\mu(\phi_n); \qquad a_n: = n r_n^2
p_n.
\]
Since we assume $d=2$, it follows from definitions (\ref{Idef})
and (\ref{J2def}) that
%
\begin{equation}
n I(\phi_n) = a_n J_2(\phi_n),\qquad n
\in\N. \label{0731b}
\end{equation}
In this section we assume
$r_n = n^{-\Omega(1)}$, so in particular
$r_n = o(1)$.

Set $N_0(n):= N_0 (G_{\phi_n}(\Po_n))$.
By the Mecke formula, $\E N_0(G_{\phi_n}(\Po_n)) =I_n(\phi_n)$,
where we set
$I_n(\phi):= n \int_{\Gamma} \exp
 (- n \int_\Gamma\phi(y-x) \,dy  ) \,dx $, so
$I_n(\phi_n) $
is
the left-hand side of (\ref{prcond1}).

Given $\eps>0 $,
truncate $\phi_n$ by setting $\tilde{\phi}_n(x):=
\phi_n(x) {\mathbf1}_{[0,r_n^{1-\eps}]} (|x|)$ for $x \in\R^2$.
Couple $G_{\phi_n}(\Po_n)$ and
$G_{\tilde{\phi}_n}(\Po_n)$
as in the proof of Lemma~\ref{poislem1}.
Let
$\tilde{N}_0(n):= N_0 (G_{\tilde{\phi}_n}(\Po_n))$.
Let
$N_0^\mathrm{ int}:= N_0^\mathrm{ int}(n)$ denote the number of isolated
vertices of
$G_{\tilde{\phi}_n} (\Po_n)$ lying in $[r^{1-\eps}_n,1-r^{1-\eps}_n]^2$.
Let
$N_0^\mathrm{ side}:= N_0^\mathrm{ side}(n)$ denote the number of isolated
vertices of
$G_{\tilde{\phi}_n} (\Po_n)$ lying within Euclidean distance
$r^{1-\eps}_n$ of precisely one edge of $\Gamma$. Let
$N_0^\mathrm{ cor}:= N_0^\mathrm{ cor}(n)$ denote the number of isolated
vertices of
$G_{\tilde{\phi}_n} (\Po_n)$ lying within $\ell_\infty$ distance
$r^{1-\eps}_n$ of one of the corners of $\Gamma$. Then
$\tilde{N}_0 (n) = N_0^\mathrm{ int}
+ N_0^\mathrm{ side}
+ N_0^\mathrm{ cor}$ (with probability 1), so
\[
I_n(\tilde{\phi}_n) = \E N_0^\mathrm{ int}
+ \E N_0^\mathrm{ side} + \E N_0^\mathrm{ cor}.
\]
Also, if $r_n= n^{-\Omega(1)}$, then
%
\begin{equation}
0 \leq\E\tilde{N}_0(n) - \E N_0(n) \leq n^2
\phi_n\bigl(r_n^{1-\eps}\bigr) \leq3 n^2
\exp\bigl( -\eta r_n^{-\eps\eta}\bigr) \to0 \label{1022a}
\end{equation}
and
%
\begin{eqnarray}\label{I2I}
n\bigl( I(\phi_n) - I (\tilde{\phi}_n)\bigr) &=&
nr_n^2 \int_{\{x\dvtx |x| \geq r_n^{-\eps}\}}
\phi_n(r_n x) \,dx
\nonumber
\\[-8pt]
\\[-8pt]
\nonumber
&\leq&3 nr_n^2 \int_{\{x\dvtx |x|> r_n^{-\eps}\}}
\eta^{-1} \exp\bigl(- \eta |x|\bigr) \,dx \to0.
\end{eqnarray}

As with (\ref{fromThetacond}), a necessary condition
for (\ref{prcond1}) is that
%
\begin{equation}
n p_n r_n^2 = \Theta(\log n). \label{nfromThetacond}
\end{equation}
Recall from (\ref{J1def}) that
$J_1(\phi_n):= J_1(\phi_n,\eta):=
p_n^{-1} \int_{0}^\infty\phi_n ((r_n t,0))
\,dt$.

\begin{lemm}
\label{lemside}
Suppose
(\ref{nfromThetacond}) holds,
and
$r_n = n^{-\Omega(1)}$
as $n \to\infty$.
Then
provided $\eps>0 $ is chosen sufficiently small (but fixed),
as $n \to\infty$ we have
%
\begin{equation}
\E N_0^\mathrm{ side} \sim \frac{2}{J_1(\phi_n)} \biggl(
\frac{ n }{a_n p_n } \biggr)^{1/2} e^{- n I(\phi_n)/2} \label{0624b}
\end{equation}
and
%
\begin{equation}
\E N_0^\mathrm{ cor} \sim\frac{4 e^{-n I(\phi_n)/4}}{a_n p_n J_1(\phi_n)^2}. \label{0702c}
\end{equation}
\end{lemm}

\begin{pf}
For $u >0$, let
\[
f_n(u):= p_n^{-1} \int_{[0,\infty) \times[0,u] }
\tilde{\phi}_n (r_n x) \,dx.
\]
Then we claim that
for $\theta_n = a_n$ or $\theta_n = 2 a_n$,
%
\begin{equation}
\int_0^{r_n^{-\eps}} \exp\bigl( - \theta_n
f_n (u)\bigr) \,du \sim1/\bigl(\theta_n J_1(
\phi_n)\bigr)\qquad \mbox{as } n \to\infty. \label{0731a}
\end{equation}
To see this, note first that $J_1(\tilde{\phi}_n) \sim J_1(\phi_n)$
as $n \to\infty$, by (\ref{Jbds}).
Also, since $\tilde{\phi}_n(x)$ is nonincreasing
in $|x|$ (because $\phi_n \in\Phi_{2,\eta} \cap\Psi_2$) we have
%
\begin{equation}
f_n(u) \leq u J_1(\tilde{\phi}_n),
\label{0905b}
\end{equation}
so that using (\ref{Jbds}) we have
%
\begin{eqnarray}\label{0926a}
&&\int_0^{r_n^{-\eps}} \exp\bigl( - \theta_n
f_n (u)\bigr) \,du \nonumber\\
&&\qquad\geq\int_0^{r_n^{-\eps}}
\exp\bigl(-\theta_n u J_1(\tilde{\phi}_n)
\bigr) \,du
\\
&&\qquad= \bigl(\theta_n J_1(\tilde{\phi}_n)
\bigr)^{-1} \int_0^{\theta_n
J_1(\tilde{\phi}_n) r_n^{-\eps}}
e^{-t}\,dt \sim \bigl( \theta_n J_1(\tilde{
\phi}_n)\bigr)^{-1}. \nonumber
\end{eqnarray}
Also given $\delta>0$, for $(s,t) \in[0,\infty) \times
(0, \delta r_n )$, we have
$\phi_n((s,t)) \geq\phi_n((s + \delta r_n,0))$, and hence
%
\begin{eqnarray}\label{1008a}
&&\int_0^\delta\exp\bigl( - \theta_n
f_n (u)\bigr) \,du\nonumber\\
&&\qquad \leq\int_0^\delta\exp
\biggl( - \theta_n u p_n^{-1} \int
_0^\infty \tilde{\phi}_n\bigl(
\bigl(r_n (s + \delta),0\bigr)\bigr) \,ds \biggr) \,du
\\
&&\qquad\leq\int_0^\delta\exp\bigl(- \theta_n
u \bigl(J_1(\tilde{\phi}_n) -\delta \bigr) \bigr) \,du \sim
\bigl(\theta_n\bigl(J_1(\tilde{\phi}_n) -
\delta\bigr) \bigr)^{-1}, \nonumber
\end{eqnarray}
and provided
$\delta\leq1/2$, we also have for $u \geq\delta$ that
%
\begin{equation}
f_n(u) \geq f_n(\delta) \geq p_n^{-1}
\int_{[0,1/2] \times[0,\delta] } \phi_n(r_n x) \,dx \geq
\delta\eta/2, \label{0926d}
\end{equation}
so that
%
\begin{equation}
\int_\delta^1 \exp\bigl( - \theta_n
f_n (u)\bigr) \,du \leq \exp(- \delta\eta\theta_n /2) = o
\bigl(\theta_n^{-1}\bigr). \label{0926b}
\end{equation}
For $u \geq1$ we have
$f_n(u) \geq f_n(1/2) \geq\eta/4$, and for $n$ large enough
$r_n^{-2} \leq n$ by (\ref{nfromThetacond}),
so
\[
\int_1^{r_n^{-\eps}} e^{-\theta_n f_n(u) } \,du \leq
r_n^{-\eps} \exp(- \eta\theta_n /4) \leq
n^{\eps/2} \exp( - \eta\theta_n /4).
\]
%
Provided $\eps$ is small enough, using (\ref{nfromThetacond}) again
we have
that the last expression is less than $\exp(- \eta\theta_n/8)$,
which is $o(\theta_n^{-1})$.
Combining this with (\ref{0926a}), (\ref{1008a}) and (\ref{0926b})
and using the fact that $\delta$ can be
arbitrarily small, we have (\ref{0731a}).

Since $\tilde{\phi}_n$ has range $r_n^{1-\eps}$
we have
\[
\E N_0^\mathrm{ side} = \bigl(4 +o(1)\bigr) n \exp\bigl(- n I(
\tilde{\phi}_n)/2\bigr) \int_0^{r_n^{-\eps}}
\exp\bigl( - 2 n r_n^2 p_n f_n(u)
\bigr) r_n \,du.
\]
By (\ref{I2I}) and
(\ref{0731a}) we obtain
\[
\E N_0^\mathrm{ side} \sim \frac{4 n r_n e^{-n I(\phi_n) /2}}{2 J_1(\phi_n) n r_n^2 p_n } = \frac{2}{J_1(\phi_n)}
\biggl( \frac{ n }{a_n p_n } \biggr)^{1/2} e^{-n I(\phi_n)/2}.
\]

Now consider $\E N_0^\mathrm{ cor}$.
For $u,v>0$, set
\[
g_n(u,v):= p_n^{-1} \int_{[0,u]\times[0,v]}
\tilde{\phi}_n\bigl(r_n\bigl(x-(u,v)\bigr) \bigr) \,dx.
\]
Then since $\tilde{\phi}_n$ has range $r_n^{1-\eps}$,
%
\begin{equation}
\E N_0^\mathrm{ cor} = \bigl(1+o(1)\bigr) 4 r_n^2
n e^{- n I(\tilde{\phi}_n) /4} \tilde{I}_n \label{0702a}
\end{equation}
with
\[
\tilde{I}_n:= \int_0^{r_n^{-\eps}} \int
_0^{r_n^{-\eps}} \exp\bigl(-n p_n
r_n^2 \bigl[ f_n (u) + f_n(v) +
g_n(u,v) \bigr] \bigr) \,du \,dv.
\]
For $u,v \geq0$ we have
$0 \leq g_n(u,v) \leq uv $. Hence
by (\ref{0905b})
we have
\begin{eqnarray*}
\tilde{I}_{n} &\geq&\int_0^{r_n^{-\eps}} \int
_0^{r_n^{-\eps}} \exp\bigl(-a_n \bigl(u
J_1(\tilde{\phi}_n) + v J_1(\tilde{
\phi}_n) + uv \bigr) \bigr) \,du \,dv
\\
&\sim& \int_0^{r_n^{-\eps}} \biggl( \frac{
e^{-a_n v J_1(\tilde{\phi}_n)} }{a_n (J_1(\tilde{\phi}_n) +v)}
\biggr) \,dv \sim\bigl(a_n J_1(\tilde{\phi}_n)
\bigr)^{-2}.
\end{eqnarray*}
On the other hand, given $\delta\in(0,\eta)$, similarly
to
(\ref{1008a}), the contribution to
$\tilde{I}_n$ from $\max(u,v) \leq\delta$ is bounded
above by
\[
\int_0^\delta\int_0^\delta
\exp \bigl(- a_n \bigl[u \bigl(J_1(\phi_n) -
\delta\bigr) + v\bigl(J_1(\phi_n ) -\delta\bigr) \bigr]
\bigr) \,du \,dv \sim\bigl(a_n \bigl(J_1(\phi_n) -
\delta\bigr)\bigr)^{-2},
\]
while by (\ref{0926d}) the contribution
to $\tilde{I}_n$ from $1 \geq\max(u,v) > \delta$ is bounded
above by $\exp(- a_n \eta\delta/2)$,
which is $o(a_n^{-2})$ by (\ref{nfromThetacond}),
and the contribution
to $\tilde{I}_n$ from $ \max(u,v) > 1$ is bounded
above by $\exp(- a_n \eta/4)r_n^{-2 \eps}$,
and hence [using (\ref{nfromThetacond})], by
$n^{\eps} \exp(- a_n \eta/4 )$,
which is $o(a_n^{-2})$ provided
$\eps$ is taken sufficiently small.
Therefore
we have $\tilde{I}_n \sim
(a_n J_1(\phi_n) )^{-2}$.
Then by (\ref{0702a})
we get (\ref{0702c}).
\end{pf}

\begin{lemm}
\label{lem1008}
Fix $\eps\in(0,1)$.
Suppose $r_n = n^{-\Omega(1)}$. Then $\E N_0^\mathrm{ int} \sim
n e^{- n I(\phi_n)}$ as $n \to\infty$.
\end{lemm}
\begin{pf}
The result follows from (\ref{I2I}).
\end{pf}

\begin{prop}
\label{lemrp1}
Suppose $d=2$.
Let $\alpha\in(0,\infty)$.
Suppose for some $\eta\in(0,1]$ that
$\phi_n \in\Phi_{2,\eta} \cap\Psi_2$ for all $n$, and
$p_n = \omega(1/\log n) $ as $n \to\infty$.
Then (\ref{prcond1}) holds if
%
\begin{equation}
n I(\phi_n) - \log n \to- \log\alpha. \label{rladef}
\end{equation}
\end{prop}

\begin{pf}
Assume (\ref{rladef}) holds, which implies {a fortiori}
that
(\ref{nfromThetacond}) also holds, so
in particular $r_n^2= O((\log n)^2/n)$. Then
by Lemma~\ref{lem1008} and (\ref{rladef}) we have
$
\E N_0^\mathrm{ int} \to\alpha$.

Using
(\ref{nfromThetacond}), (\ref{rladef}) and
Lemma~\ref{lemside},
we obtain (for a sufficiently small choice of $\eps$)
that
$\E N_0^\mathrm{ side}
= O( (p_n \log n)^{-1/2} ) $,
which tends to zero by the assumption on
$p_n $.
Similarly,
by (\ref{0702c}) and (\ref{rladef}),
$
\E N_0^\mathrm{ cor} = O(n^{-1/4} /(p_n \log n)) = o(1)$.
Applying (\ref{1022a})
completes the proof.
\end{pf}

When $p_n = O(1/\log n)$,
boundary
effects become important in the asymptotics for
the mean number
of isolated points.

\begin{prop}
\label{lemr2}
Suppose $d=2$, and for some $\eta\in(0,1]$ we
have $\phi_n \in\Phi_{2,\eta} \cap\Psi_2$ for all $n$.
Suppose $p_n = o(1/\log n)$
and also $p_n = \omega((\log n)^{-1}n^{-1/3} )$,
as $n \to\infty$. Fix $\alpha\in(0,\infty)$, and assume
%
\begin{equation}
n I(\phi_n) = \log \biggl( \frac{4 J_2(\phi_n)
}{\alpha^2J_1(\phi_n)^2} \biggr) + \log \biggl(
\frac{n}{p_n} \biggr) - 
\log\log \biggl( \frac{n}{p_n} \biggr)
+ o(1). \label{ala}
\end{equation}
Then (\ref{prcond1}) holds.
\end{prop}

\begin{pf}
Under the assumptions given,
using Lemma~\ref{lem1008} we have
%
\begin{equation}
\E N_0^\mathrm{ int} = \bigl(1 +o(1)\bigr) n e^{-n I(\phi_n)} = O
\biggl( p_n \log \biggl( \frac{ n }{p_n} \biggr) \biggr) \to0.
\label{fromala}
\end{equation}
Also by (\ref{0624b}), (\ref{ala}) and (\ref{0731b}),
we have
%
\begin{equation}
\E N_0^\mathrm{ side} \sim\alpha \biggl( \frac{ n}{J_2(\phi_n) a_n p_n }
\biggr)^{1/2} \biggl( \frac{p_n}{n} \biggr)^{1/2} \bigl(
\log(n /p_n) \bigr)^{1/2} \to\alpha. \label{0710b}
\end{equation}

Using
(\ref{0624b}) again along with (\ref{Jbds}), we obtain
that $e^{-n I(\phi_n)/4} = \Theta((a_n p_n/ \break n)^{1/4})$ so that
(\ref{0702c}) yields
$
\E N_0^\mathrm{ cor}
= O (  ((a_n p_n)^3 n  )^{-1/4})$,
and by (\ref{nfromThetacond})
[which follows from (\ref{ala})]
and the assumption
$p_n = \omega(n^{-1/3} (\log n)^{-1})$, this shows that
$\E N_0^\mathrm{ cor} \to0$.
Combined with (\ref{fromala}), (\ref{0710b}) and (\ref{1022a}),
we have the result.
\end{pf}

Consider the intermediate case with $p_n = \Theta(1/ \log n)$.

\begin{theo}
\label{lemrp5}
Let $\alpha\in(0,\infty), \eta\in(0,1]$.
Suppose that $\phi_n \in\Phi_{2,\eta} \cap\Psi_2$ for all $n$,
$p_n = \Theta(1/ \log n)$ and
$n I(\phi_n) = \log n -2 \log\gamma_n + o(1)$,
where $\gamma_n$ denotes the solution in $(0,\infty)$ to
%
\begin{equation}
\gamma_n^2 + 2 \gamma_n
\bigl(J_2(\phi_n)^{1/2}/J_1(
\phi_n)\bigr) (p_n \log n)^{-1/2} = \alpha.
\label{0703a}
\end{equation}
Then (\ref{prcond1}) holds.
\end{theo}

\begin{pf}
By (\ref{Jbds}) and the assumption on $p_n$,
$\limsup_{n \to\infty}(\gamma_n) < \infty$
and
$\liminf_{n \to\infty}(\gamma_n) > 0$.
By Lemma~\ref{lem1008},
\[
\E N_0^\mathrm{ int} = \bigl(1+o(1)\bigr) n e^{-n I(\phi_n)} =
\bigl(1+o(1)\bigr) \gamma_n^{2},
\]
while by (\ref{0624b}) and (\ref{0731b}),
\[
\E N_0^\mathrm{ side} \sim \frac{2}{J_1(\phi_n)} \biggl(
\frac{ n}{a_n p_n} \biggr)^{1/2} \gamma_n n^{-1/2}
\sim\frac{ 2 J_2(\phi_n)^{1/2} \gamma_n}{ J_1(\phi_n) ( p_n \log
n)^{1/2} }.
\]
Also by (\ref{0702c}), $\E N_0^\mathrm{ cor} = O(n^{-1/4}/(p_n \log n))= o(1)$.
Combining
these results and using
(\ref{0703a}) and (\ref{1022a})
we have (\ref{prcond1}).
\end{pf}

In the case $p_n = o((\log n)^{-1}n^{-1/3})$ the main contribution
to $\E N_0 $ comes from near the corners of $\Gamma$.

\begin{prop}
\label{lemrp6}
Suppose $d =2$.
Let $\alpha\in(0,\infty), \eta\in(0,1]$, and
suppose $(\phi_n)_{n >0}$ are such that
$\phi_n \in\Phi_{2,\eta} \cap\Psi_2$ for all $n$ and
$p_n = o((\log n)^{-1}n^{-1/3} )$ and
%
\begin{eqnarray}\label{0703b}
n I(\phi_n) &=& 4 \bigl( \log(1/p_n) 
- \log
\log(1/p_n) + \log\bigl(J_2(\phi_n)/ \bigl(
\alpha J_1(\phi_n)^2\bigr) \bigr) \bigr)
\nonumber
\\[-8pt]
\\[-8pt]
\nonumber
&&{} +o(1)
\end{eqnarray}
as $n \to\infty$.
Assume also that
$r_n = n^{-\Omega(1)}$.
Then (\ref{prcond1}) holds.
\end{prop}

\begin{pf}
Note that $p_n = \Omega((\log n)/n)$ since
otherwise (\ref{0703b}) cannot be satisfied by bounded $r_n$.
Then $\log(1/ p_n) = \Theta(\log n)$ and (\ref{nfromThetacond}) holds.
By (\ref{0702c}) and (\ref{0731b}),
%
\begin{equation}
\E N_0^\mathrm{ cor} \sim\frac{4 p_n \log(1/p_n)
\alpha J_1(\phi_n)^2/ J_2(\phi_n) }{J_1(\phi_n)^2 a_n p_n} \to\alpha.
\label{0703d}
\end{equation}
Also
$e^{-n I(\phi_n)/4} = \Theta(p_n \log(1/p_n))= \Theta(p_n \log n)$.
Therefore by (\ref{0624b}) and (\ref{nfromThetacond}),
we obtain that
$
\E N_0^\mathrm{ side}
= O  ( n^{1/2} (p_n \log n)^{3/2}  )$,
which tends to zero since
we assume
$p_n = o(n^{-1/3} (\log n)^{-1})$.

Finally, since $p_n = \Omega((\log n)/n)$, using
Lemma~\ref{lem1008} and (\ref{0703b}) we have
\[
\E N_0^\mathrm{ int} = O\bigl(n e^{-n I(\phi_n)}\bigr) = O\bigl(
n p_n^4 (\log1/p_n)^4\bigr) = O
\bigl(n^{-2}\bigr),
\]
so
$\E N_0^\mathrm{ int} \to0$, and (\ref{prcond0}) then follows by
(\ref{1022a}).
\end{pf}

\begin{pf*}{Proof of Theorem~\ref{thconds}}
The proof follows immediately from Propositions~\ref{lemrp1}, \ref
{lemr2} and \ref{lemrp6}.
\end{pf*}

Our final result deals with the intermediate case with
$p_n= \break \Theta(n^{-1/3} (\log n)^{-1})$.


\begin{theo}
\label{lemrp4}
Let $\alpha\in(0,\infty)$, and
suppose $(\phi_n)_{n >0}$ are such that
$p_n = \Theta(n^{-1/3} (\log n)^{-1})$ and
%
\begin{eqnarray}\label{0703c}
n I(\phi_n)& =& 4 \bigl( \log(1/p_n) 
- \log
\log(1/p_n) + \log\bigl(J_2(\phi_n)/\bigl(
\beta_n J_1(\phi_n)^2\bigr)\bigr)
\bigr)
\nonumber
\\[-8pt]
\\[-8pt]
\nonumber
&&{} +o(1)
\end{eqnarray}
as $n \to\infty$, with $\beta_n$ denoting the solution in
$(0,\infty)$ to
%
\begin{equation}
\bigl(3 J_2(\phi_n)\bigr)^{-3/2}
J_1(\phi_n)^{3} \bigl(n^{1/3}
p_n \log n\bigr)^{3/2} \beta_n^2 +
\beta_n = \alpha. \label{0711a}
\end{equation}
Then (\ref{prcond1}) holds.
\end{theo}

\begin{pf}
Note that (\ref{0703c}) is the same as (\ref{0703b}) but
with $\alpha$ replaced by $\beta_n$.
As with (\ref{0703d}) we have $\E N_0^\mathrm{ cor} = \beta_n +o(1)$.
Then by (\ref{0624b}) and (\ref{0703c})
we have
\[
\E N_0^\mathrm{ side} \sim \bigl(2/J_1(
\phi_n)\bigr) (n/a_n)^{1/2} p_n^{3/2}
(\log1/p_n)^2 \beta_n^2
J_1(\phi_n)^4 J_2(
\phi_n)^{-2}.
\]
By (\ref{0731b}) and (\ref{0703c}),
$
a_n = n I(\phi_n)/J_2(\phi_n) \sim(4/J_2(\phi_n))
\log(1/p_n)$,
and our assumption on $p_n$ implies $\log1/p_n \sim(1/3) \log n$,
so that
\[
\E N_0^\mathrm{ side} \sim \beta_n^2
J_1(\phi_n)^{3} J_2(
\phi_n)^{-3/2} n^{1/2} p_n^{3/2}
\bigl( (\log n)/3\bigr)^{3/2} .
\]
Hence by (\ref{0711a}), $\E[N_0^\mathrm{ side} + N_0^\mathrm{ cor}] \to
\alpha$.
Also by Lemma~\ref{lem1008},
$\E N_0^\mathrm{ int} = \break  O(n e^{-n I(\phi_n)}) = O(n p_n^4(\log1/p_n)^4)$,
which tends to zero, and (\ref{prcond1}) follows by~(\ref{1022a}).
\end{pf}






\printaddresses

\begin{thebibliography}{18}


\bibitem{BBKMW}
\begin{barticle}[mr]
\bauthor{\bsnm{Balogh},~\bfnm{J{\'o}zsef}\binits{J.}},
\bauthor{\bsnm{Bollob{\'a}s},~\bfnm{B{\'e}la}\binits{B.}},
\bauthor{\bsnm{Krivelevich},~\bfnm{Michael}\binits{M.}},
\bauthor{\bsnm{M{\"u}ller},~\bfnm{Tobias}\binits{T.}} \AND
\bauthor{\bsnm{Walters},~\bfnm{Mark}\binits{M.}}
(\byear{2011}).
\btitle{Hamilton cycles in random geometric graphs}.
\bjournal{Ann. Appl. Probab.}
\bvolume{21}
\bpages{1053--1072}.
\bid{doi={10.1214/10-AAP718}, issn={1050-5164}, mr={2830612}}
\end{barticle}
%

\bptok{imsref}%
\endbibitem

\bibitem{Boll}
\begin{bbook}[mr]
\bauthor{\bsnm{Bollob{\'a}s},~\bfnm{B{\'e}la}\binits{B.}}
(\byear{2001}).
\btitle{Random Graphs},
\bedition{2nd} ed.
\bseries{Cambridge Studies in Advanced Mathematics}
\bvolume{73}.
\bpublisher{Cambridge Univ. Press},
\blocation{Cambridge}.
\bid{doi={10.1017/CBO9780511814068}, mr={1864966}}
\end{bbook}
%

\bptok{imsref}%
\endbibitem

\bibitem{BDFL}
\begin{barticle}[mr]
\bauthor{\bsnm{Broutin},~\bfnm{Nicolas}\binits{N.}},
\bauthor{\bsnm{Devroye},~\bfnm{Luc}\binits{L.}},
\bauthor{\bsnm{Fraiman},~\bfnm{Nicolas}\binits{N.}} \AND
\bauthor{\bsnm{Lugosi},~\bfnm{G{\'a}bor}\binits{G.}}
(\byear{2014}).
\btitle{Connectivity threshold of {b}luetooth graphs}.
\bjournal{Random Structures Algorithms}
\bvolume{44}
\bpages{45--66}.
\bid{doi={10.1002/rsa.20459}, issn={1042-9832}, mr={3143590}}
\end{barticle}
%

\bptok{imsref}%
\endbibitem

\bibitem{CDG}
\begin{barticle}[mr]
\bauthor{\bsnm{Coon},~\bfnm{Justin}\binits{J.}},
\bauthor{\bsnm{Dettmann},~\bfnm{Carl~P.}\binits{C.~P.}} \AND
\bauthor{\bsnm{Georgiou},~\bfnm{Orestis}\binits{O.}}
(\byear{2012}).
\btitle{Full connectivity: Corners, edges and faces}.
\bjournal{J. Stat. Phys.}
\bvolume{147}
\bpages{758--778}.
\bid{doi={10.1007/s10955-012-0493-y}, issn={0022-4715}, mr={2930579}}
\end{barticle}
%

\bptok{imsref}%
\endbibitem

\bibitem{DPS}
\begin{barticle}[mr]
\bauthor{\bsnm{Diaz},~\bfnm{Josep}\binits{J.}},
\bauthor{\bsnm{Petit},~\bfnm{Jordi}\binits{J.}} \AND
\bauthor{\bsnm{Serna},~\bfnm{Maria}\binits{M.}}
(\byear{2000}).
\btitle{Faulty random geometric networks}.
\bjournal{Parallel Process. Lett.}
\bvolume{10}
\bpages{343--357}.
\bid{doi={10.1142/S0129626400000329}, issn={0129-6264}, mr={1841124}}
\end{barticle}
%

\bptok{imsref}%
\endbibitem

\bibitem{ER}
\begin{barticle}[mr]
\bauthor{\bsnm{Erd{\H{o}}s},~\bfnm{P.}\binits{P.}} \AND
\bauthor{\bsnm{R{\'e}nyi},~\bfnm{A.}\binits{A.}}
(\byear{1959}).
\btitle{On random graphs. {I}}.
\bjournal{Publ. Math. Debrecen}
\bvolume{6}
\bpages{290--297}.
\bid{issn={0033-3883}, mr={0120167}}
\end{barticle}
%

\bptok{imsref}%
\endbibitem

\bibitem{GK}
\begin{bincollection}[mr]
\bauthor{\bsnm{Gupta},~\bfnm{Piyush}\binits{P.}} \AND
\bauthor{\bsnm{Kumar},~\bfnm{P.~R.}\binits{P.~R.}}
(\byear{1999}).
\btitle{Critical power for asymptotic connectivity in wireless networks}.
In \bbooktitle{Stochastic Analysis, Control, Optimization and Applications}
\bpages{547--566}.
\bpublisher{Birkh\"auser},
\blocation{Boston, MA}.
\bid{mr={1702981}}
\end{bincollection}
%

\bptok{imsref}%
\endbibitem

\bibitem{GK2}
\begin{barticle}[mr]
\bauthor{\bsnm{Gupta},~\bfnm{Piyush}\binits{P.}} \AND
\bauthor{\bsnm{Kumar},~\bfnm{P.~R.}\binits{P.~R.}}
(\byear{2000}).
\btitle{The capacity of wireless networks}.
\bjournal{IEEE Trans. Inform. Theory}
\bvolume{46}
\bpages{388--404}.
\bid{doi={10.1109/18.825799}, issn={0018-9448}, mr={1748976}}
\end{barticle}
%

\bptok{imsref}%
\endbibitem

\bibitem{KGM}
\begin{bmisc}[auto:parserefs-M02]
\bauthor{\bsnm{Krishnan},~\bfnm{B.~S.}\binits{B.~S.}},
\bauthor{\bsnm{Ganesh},~\bfnm{A.}\binits{A.}} \AND
\bauthor{\bsnm{Manjunath},~\bfnm{D.}\binits{D.}}
(\byear{2013}).
\bhowpublished{%
On connectivity thresholds in superposition of random key graphs
on random geometric graphs.
In \textit{Information Theory Proceedings (ISIT),
2013 IEEE International Symposium on 7--12 July 2013}
2389--2393. IEEE, New York}.
\end{bmisc}
%

\bptok{imsref}%
\endbibitem

\bibitem{MaoAnderson}
\begin{barticle}[auto:parserefs-M02]
\bauthor{\bsnm{Mao},~\bfnm{G.}\binits{G.}} \AND
\bauthor{\bsnm{Anderson},~\bfnm{B.~D.~O.}\binits{B.~D.~O.}}
(\byear{2012}).
\btitle{Towards a better understanding of large-scale network models}.
\bjournal{IEEE/ACM Transactions on Networking}
\bvolume{20}
\bpages{408--421}.
\end{barticle}
%

\bptok{imsref}%
\endbibitem

\bibitem{MR}
\begin{bbook}[mr]
\bauthor{\bsnm{Meester},~\bfnm{Ronald}\binits{R.}} \AND
\bauthor{\bsnm{Roy},~\bfnm{Rahul}\binits{R.}}
(\byear{1996}).
\btitle{Continuum Percolation}.
\bseries{Cambridge Tracts in Mathematics}
\bvolume{119}.
\bpublisher{Cambridge Univ. Press},
\blocation{Cambridge}.
\bid{doi={10.1017/CBO9780511895357}, mr={1409145}}
\end{bbook}
%

\bptok{imsref}%
\endbibitem

\bibitem{Penbk}
\begin{bbook}[mr]
\bauthor{\bsnm{Penrose},~\bfnm{Mathew}\binits{M.}}
(\byear{2003}).
\btitle{Random Geometric Graphs}.
\bseries{Oxford Studies in Probability}
\bvolume{5}.
\bpublisher{Oxford Univ. Press},
\blocation{Oxford}.
\bid{doi={10.1093/acprof:oso/9780198506263.001.0001}, mr={1986198}}
\end{bbook}
%

\bptok{imsref}%
\endbibitem

\bibitem{conperc}
\begin{barticle}[mr]
\bauthor{\bsnm{Penrose},~\bfnm{Mathew~D.}\binits{M.~D.}}
(\byear{1991}).
\btitle{On a continuum percolation model}.
\bjournal{Adv. in Appl. Probab.}
\bvolume{23}
\bpages{536--556}.
\bid{doi={10.2307/1427621}, issn={0001-8678}, mr={1122874}}
\end{barticle}
%

\bptok{imsref}%
\endbibitem

\bibitem{MSTLong}
\begin{barticle}[mr]
\bauthor{\bsnm{Penrose},~\bfnm{Mathew~D.}\binits{M.~D.}}
(\byear{1997}).
\btitle{The longest edge of the random minimal spanning tree}.
\bjournal{Ann. Appl. Probab.}
\bvolume{7}
\bpages{340--361}.
\bid{doi={10.1214/aoap/1034625335}, issn={1050-5164}, mr={1442317}}
\end{barticle}
%

\bptok{imsref}%
\endbibitem

\bibitem{kcon}
\begin{barticle}[mr]
\bauthor{\bsnm{Penrose},~\bfnm{Mathew~D.}\binits{M.~D.}}
(\byear{1999}).
\btitle{On {$k$}-connectivity for a geometric random graph}.
\bjournal{Random Structures Algorithms}
\bvolume{15}
\bpages{145--164}.
\bid{doi={10.1002/(SICI)1098-2418(199909)15:2<145::AID-RSA2>3.0.CO;2-G}, issn={1042-9832}, mr={1704341}}
\end{barticle}
%

\bptok{imsref}%
\endbibitem

\bibitem{Tse}
\begin{bbook}[auto:parserefs-M02]
\bauthor{\bsnm{Tse},~\bfnm{D.}\binits{D.}} \AND
\bauthor{\bsnm{Viswanath},~\bfnm{P.}\binits{P.}}
(\byear{2005}).
\btitle{Fundamentals of Wireless Communication}.
\bpublisher{Cambridge Univ. Press},
\blocation{Cambridge}.
\end{bbook}
%

\bptok{imsref}%
\endbibitem

\bibitem{Yagan}
\begin{barticle}[mr]
\bauthor{\bsnm{Ya{\u{g}}an},~\bfnm{Osman}\binits{O.}}
(\byear{2012}).
\btitle{Performance of the {E}schenauer--{G}ligor key distribution scheme under an ON/{OFF} channel}.
\bjournal{IEEE Trans. Inform. Theory}
\bvolume{58}
\bpages{3821--3835}.
\bid{doi={10.1109/TIT.2012.2189353}, issn={0018-9448}, mr={2924403}}
\end{barticle}
%

\bptok{imsref}%
\endbibitem

\bibitem{YWLH}
\begin{bmisc}[auto:parserefs-M02]
\bauthor{\bsnm{Yi},~\bfnm{C.-W.}\binits{C.-W.}},
\bauthor{\bsnm{Wan},~\bfnm{P.-J.}\binits{P.-J.}},
\bauthor{\bsnm{Lin},~\bfnm{K.-W.}\binits{K.-W.}} \AND
\bauthor{\bsnm{Huang},~\bfnm{C.-H.}\binits{C.-H.}}
(\byear{2006})
\bhowpublished{%
Asymptotic distribution of the number of isolated nodes in wireless ad hoc networks with unreliable nodes and links.
In \textit{Global Telecommunications Conference
2006, GLOBECOM'06}. IEEE, New York}.
\end{bmisc}
%

\bptok{imsref}%
\endbibitem
\end{thebibliography}
\end{document}